\documentclass[hidelinks,onefignum,onetabnum]{siamart220329}

\usepackage[english]{babel} 
\usepackage{amssymb,amsfonts,txfonts}
\usepackage{graphicx,xcolor} 
\usepackage{epstopdf}
\usepackage[]{algorithmic}

\ifpdf
  \DeclareGraphicsExtensions{.eps,.pdf,.png,.jpg}
\else
  \DeclareGraphicsExtensions{.eps}
\fi

\usepackage{comment} 

\newcommand{\IN}{\mathbb{N}}
\newcommand{\IR}{\mathbb{R}}
\newcommand{\IRn}{\mathbb{R}^n}
\newcommand{\cR}{\mathcal{R}}
\newcommand{\cP}{\mathcal{P}}
\newcommand{\cD}{\mathcal{D}}

\newcommand{\cL}{\mathcal{L}}

\newcommand{\cS}{\mathcal{S}}

\newcommand{\cO}{\mathcal{O}}
\DeclareMathOperator*{\argmin}{\textrm{argmin}}
\DeclareMathOperator*{\argmax}{\textrm{argmax}}
\newcommand{\prox}{\textrm{prox}}

\newcommand{\dom}{\textrm{dom}\,}
\newcommand{\ri}{\textrm{ri}}
\newcommand{\MGOPT}{\texttt{MGOPT} }
\newcommand{\MGProx}{\texttt{MGProx} }

\def\D {\mathbf{D}}
\def\E {\mathbf{E}}
\def\F {\mathbf{F}}
\def\P {\mathbf{P}}
\def\R {\mathbf{R}}
\def\U {\mathbf{U}}
\def\c {\mathbf{c}}
\def\u {\mathbf{u}}
\def\v {\mathbf{v}}
\def\z {\mathbf{z}}
\def\bphi {\mathbf{\phi}}
\def\bpsi {\mathbf{\psi}}
\def\vec {\textrm{vec}}
\def\grad {\textbf{grad}}
\def\Hess {\textbf{Hess}}

\newsiamremark{remark}{Remark}
\newsiamremark{hypothesis}{Hypothesis}
\crefname{hypothesis}{Hypothesis}{Hypotheses}
\newsiamthm{claim}{Claim}

\newsiamthm{example}{Example}

\newlength{\commentWidth}
\headers{MGProx: A nonsmooth multigrid proximal gradient method}{A. Ang, H. De Sterck and S. Vavasis}

\title{MGProx: A nonsmooth multigrid proximal gradient method with adaptive restriction
for strongly convex optimization
    \thanks{Submitted to the editors DATE.
            \funding{
            Supported in part by a joint postdoctoral fellowship by 
            the Fields Institute for Research in Mathematical Sciences and the University of Waterloo, and in part by Discovery Grants from the Natural Sciences and Engineering Research Council (NSERC) of Canada.
            }
    }
}

\author{
Andersen Ang\thanks{Electronics and Computer Science, University of Southampton, 
United Kingdom
(\email{andersen.ang@soton.ac.uk}).
Most of the work of this paper was done when Andersen Ang was a postdoctoral fellow at the University of Waterloo.}
\and 
Hans De Sterck\thanks{Department of Applied Mathematics, University of Waterloo,
Canada
(\email{hans.desterck@uwaterloo.ca}).
                     }
\and 
Stephen Vavasis\thanks{Department of Combinatorics and Optimization, University of Waterloo, 
Canada
  (\email{vavasis@uwaterloo.ca}).
        }
}

\usepackage{amsopn}

\usepackage{soul} 
\usepackage{multirow}

\ifpdf
\hypersetup{
  pdftitle={MGProx: A multigrid proximal gradient method},
  pdfauthor={A. Ang, H. De Sterck, and  S. Vavasis}
}
\fi

\begin{document}

\maketitle

\begin{abstract}
We study the combination of proximal gradient descent with multigrid for solving a class of possibly nonsmooth strongly convex optimization problems.
We propose a multigrid proximal gradient method called MGProx, which accelerates the proximal gradient method by multigrid, based on using hierarchical information of the optimization problem.
MGProx applies a newly introduced adaptive restriction operator to simplify the Minkowski sum of subdifferentials of the nondifferentiable objective function across different levels.
We provide a theoretical characterization of MGProx.
First we show that the MGProx update operator exhibits a fixed-point property.
Next, we show that the coarse correction is a descent direction for the fine variable of the original fine level problem in the general nonsmooth case.
Lastly, under some assumptions we provide the convergence rate for the algorithm.
In the numerical tests on the Elastic Obstacle Problem, which is an example of nonsmooth convex optimization problem where multigrid method can be applied, we show that MGProx has a faster convergence speed than 
competing methods.
\end{abstract}



\section{Introduction}\label{sec:mgopt}
We study the combination of two iterative algorithms: proximal gradient descent and multigrid, to solve the following class of optimization problems
\begin{equation}\label{prob:minfg}
\argmin_x ~ F_0(x) \coloneqq f_0(x) + g_0(x).
\end{equation}
\textbf{(Assumption)} We assume $f_0 : \IRn \rightarrow \IR$ is 
$L_0$-smooth and $\mu_0$-strongly convex,
and  
$g_0 : \IRn \rightarrow \IR$ 
is proper, possibly nonsmooth, convex, 
and separable.
Regarding the non-smoothness (nondifferentiability) of $g_0$, we further assume a single point of non-differentiability for $g_0$.
We recall that a function $\zeta(x) : \IRn \rightarrow \IR$ is separable if $\zeta(x) = \sum_i \zeta_i(x_i)$, and  is $\mu$-strongly convex (and thus coercive) with $\mu>0$ if $\zeta(x) - \frac{\mu}{2}\| x \|_2^2$ is convex;
and lastly is $L$-smooth if  $\zeta$ is $C^{1,1}_L$ and $\nabla \zeta$ is $L$-Lipschitz; i.e., for all $x,y$ in $\IRn$, $\nabla \zeta(x)$ exists and
\begin{equation}
\zeta(y) \leq \zeta(x) + \big\langle \nabla \zeta(x), y-x \big\rangle + \frac{L}{2}\|y-x\|_2^2.
\label{def:Lsmooth}
\end{equation}

\paragraph{Modern models are nonsmooth}
Advancements in nonsmooth (i.e., nondifferentiable) optimization since the 60s \cite{moreau1962fonctions} enable the use of nonsmooth $g_0$ in \eqref{prob:minfg}.
The standard textbooks are \cite{rockafellar1970convex,shor1985minimization,bauschke2011convex,beck2017first}.
Then \eqref{prob:minfg} captures many models in machine learning \cite{combettes2005signal,parikh2014proximal}, where $f_0$ is a data fitting term and $g_0$ models the constraint(s) and/or regularization(s) of the application.
A popular tool for solving \eqref{prob:minfg} is the proximal gradient method \cite{passty1979ergodic,fukushima1981generalized}, to be reviewed in \cref{sec:intro:subsec:prox}.

\paragraph{Classical problems in scientific computing are smooth}
Setting $g \equiv 0$ in \eqref{prob:minfg} gives
\begin{equation}\label{prob:mgopt}
\min_{x \in \IRn} ~ f_0(x),
\end{equation}
in which this problem class of smooth strongly convex optimization subsumes many problems in scientific computing.
If problem \eqref{prob:mgopt} comes from the discretization of certain classes of partial differential equation (PDE) problems, 
 multigrid methods
\cite{fedorenko1962relaxation,brandt1977multi,hackbusch1980convergence,hackbusch1985multi,nash2000multigrid}, to be reviewed in \cref{sec:intro:subsec:mgopt}, are among the fastest known method for solving \eqref{prob:mgopt}.

\paragraph{This work: bridging smoothness and nonsmoothness}
Multigrid and nonsmooth optimization are two communities that seldom interact.
In this work we link the two fields and develop a method that can handle nonsmooth problems while enjoying the fast convergence from multigrid. 
We propose \MGProx that accelerates the proximal-gradient method by multigrid to solve Problem \eqref{prob:minfg}.
Below we review multigrid and the proximal gradient method.

%
%
\subsection{Classical multigrid and notation}\label{sec:intro:subsec:mgopt}
Multigrid dates back to the 1960s with works by Fedorenko
\cite{fedorenko1962relaxation} on solving the Poisson equation and was then further developed by 
Brandt \cite{brandt1977multi} and Hackbusch \cite{hackbusch1980convergence}.
There are many multigrid frameworks; in this work we focus on \texttt{MGOPT}: a full approximation scheme \cite{brandt1977multi} nonlinear multigrid method which was applied and extended to optimization problems by Nash \cite{nash2000multigrid}.
\MGOPT speeds up the convergence of an iterative algorithm (called smoothing or relaxation) by using a hierarchy of coarse discretizations of $f_0$:
it first constructs a series of coarse auxiliary problems of the form
\begin{equation} 
\min_{x_\ell} \,
f_\ell(x_\ell) - \langle \tau_{\ell-1 \rightarrow \ell}, \,x_\ell \rangle, 
\qquad \ell \in \{ 1,2,\dots, L \},
\label{prob:mgopt_coarse}
\end{equation}
where $\tau_{\ell-1 \rightarrow \ell}$ carries information from level $\ell-1$ to level $\ell$, 
and $f_\ell, x_\ell$ denote the function $f$  and the variable $x$, at the level $\ell$, respectively.
\MGOPT then makes use of the solution of \eqref{prob:mgopt_coarse} to solve \eqref{prob:mgopt}.
The convergence of the overall algorithm is sped up by the correction from the coarse levels and by the fact that the coarse problems are designed to be ``less expensive'' to solve than the given ones.

\paragraph{Notation}
The symbol $x_0$ (or $x$) is called the fine variable.
The symbol $x_\ell$ in \eqref{prob:mgopt_coarse} with $\ell \geq 1$ is called coarse variable.
The subscript $\ell \in \{0,1,\dots, L \}$ denotes the level.
A larger $\ell$ means a coarser level with lower resolution (fewer variables).
For the remainder of the paper, $L$ without a subscript stands for the number of levels, whereas $L_{\ell}$ denotes the smoothness parameter for the level-$\ell$ problem.
At a level $\ell$, the coarse version of the vector $x_\ell \in \IR^{n_\ell}$ is $x_{\ell+1} = \overline{\cR}(x_\ell) \coloneqq \overline{R}x_\ell$ where $\overline{R} \in \IR^{n_{\ell+1} \times n_{\ell}}$ with $n_{\ell+1} \leq n_{\ell}$ is called a restriction matrix.
Similarly, given $x_{\ell+1} \in \IR^{n_{\ell+1}}$ and a prolongation matrix $\overline{P} \in \IR^{ n_{\ell} \times n_{\ell +1}}$, 
we obtain the level-$\ell$ version of $x_{\ell+1}$ as $x_{\ell} = \overline{\cP}(x_{\ell+1}) \coloneqq \overline{P} x_{\ell+1}$.
We let $\overline{P} = c \overline{R}^\top$ given scaling factor $c>0$ and $\overline{R}$.
In multigrid, choosing $(\overline{R},\overline{P})$ depends on the application.
In this work, for the applications we consider commonly chosen $(\overline{R},\overline{P})$, see Section\,\ref{sec:exp:subsec:mgprox_eop} for our choice of 
$(\overline{R},\overline{P})$ used in the experiment and see \cite{briggs2000multigrid} for an introduction on how to select $(\overline{R},\overline{P})$.

\subsection{MGOPT}
Let $x_{\ell}^{k}$ be the level-$\ell$ variable at iteration $k$.
Algorithm~\ref{algo:mgopt} shows a 2-level ($\ell \in \{0,1\}$) \MGOPT \cite{nash2000multigrid} for solving \eqref{prob:mgopt}, with the steps in the algorithm explained as follows:
\begin{itemize}
\item (i): $\sigma: \IRn \rightarrow \IRn$ denotes an update iteration called pre-smoothing.
In this work we focus on $\sigma$ being the proximal gradient operator.

\item (ii): the restriction step.

\item (iii): the vector $\tau_{0 \rightarrow 1}^{k+1}$ carries the information at level $\ell = 0$ to level $\ell = 1$.

\item (iv): the coarse problem \eqref{prob:mgopt_coarse} is a ``smaller version'' of the original fine problem.
The function $f_1 = \overline{\cR}(f_0)$ is the coarsening of $f_0$ and the linear term $\langle \tau_{0 \rightarrow 1}^{k+1},  \xi \rangle$ links the coarse variable with the $\tau$-correction information from the fine variable.

\item (v): the updated coarse variable $x_1^{k+1}$ is used to update the fine variable $y_0^{k+1}$.

\item (vi): this step is the same as (i).
\end{itemize}
In the algorithm, $\alpha >0$ is a stepsize.
The $\tau$-correction is designed in a way that 
the iteration has a fixed-point that corresponds to a solution.

\begin{algorithm}
\caption{2-level \MGOPT \cite{nash2000multigrid} for an approximate solution of \eqref{prob:mgopt}}
\label{algo:mgopt}
\begin{algorithmic}
\STATE{Initialize $x_0^1$, $\overline{R}$ and $\overline{P}$}
\FOR{$k = 1,2,\dots$}
\STATE{(i)~~~~ $y^{k+1}_0 = \sigma(x^k_0)$ 
\hfill pre-smoothing \vspace{0.6mm}
}
\STATE{(ii)~~~ $y^{k+1}_{1}  = \overline{R} y_0^{k+1}$   
\hfill coarse variable \vspace{0.6mm}
}
\STATE{(iii)~~ $\tau^{k+1}_{0 \rightarrow 1} = \nabla f_1( y^{k+1}_1 )  - \overline{R} \nabla f_0( y^{k+1}_0 )$
\hfill create the tau vector \vspace{0.6mm}
}
\STATE{(iv)~~ $x^{k+1}_1 = \displaystyle \argmin_{\xi} \, f_1(\xi) - \langle \tau_{0\rightarrow 1}^{k+1}, \xi \rangle$
\hfill solve the coarse problem (exactly)
}
\STATE{(v)\,~~~ $z^{k+1}_0 = y^{k+1}_0  + \alpha \overline{P}\big(  x^{k+1}_1 - y_1^{k+1} \big)$
\hfill coarse correction \vspace{0.6mm}
}
\STATE{(vi)~\,~ $x^{k+1}_0 = \sigma (z^{k+1}_0)$
\hfill post-smoothing~\,
}
\ENDFOR
\end{algorithmic}
\end{algorithm}

\begin{remark}[\MGOPT has no theoretical convergence guarantee]\label{remark:NashWrong}
The proof of \cite[Theorem 1]{nash2000multigrid} on the convergence of \MGOPT requires additional assumptions.
In short the proof states the following: on solving \eqref{prob:mgopt} with an iterative algorithm $x^{k+1} \coloneqq \sigma (x^k)$ where the update map $\sigma : \IRn \rightarrow \IRn$ is assumed to be converging from any starting point $x^1$, now suppose $\rho : \IRn \rightarrow \IRn$ is some other operator with the descending property that $f_0 (\rho(x)) \leq f_0(x)$.
Then \cite[Theorem 1]{nash2000multigrid} claimed that an algorithm consisting of interlacing $\sigma$ with $\rho$ repeatedly is also convergent.
This is generally not true without further assumptions.
Here is a counterexample for \cite[Theorem 1]{nash2000multigrid}.
Consider minimizing a scalar function $f(x) =  x^2 \exp(-x^2)$. 
\begin{itemize}
    \item This $f$ has a unique global minimum at $x=0$, two global maxima at $x = \pm1$.
    \item We decrease $f$ by $\rho$ being the gradient descent step.
    \item $f$ is differentiable, and its slope is $f^\prime (x) = 2 \exp(-x^2) x (1-x^2)$ with a Lipschitz constant about $0.58$, thus we can pick $\alpha = 1 < 1/0.58$ for the gradient stepsize.
    \item If we initialize at $x^0=2$, we have $x^1= \rho(x^0) = 2 + 12e^{-4}$.
    \item Take the operator $\sigma : x \mapsto \beta x$ with  $\beta = 1/(1+6e^{-4}) = x^0/x^1$.
    For any $x$, the sequence $x, \sigma(x), \sigma(\sigma(x)),...$ converges to $x^*=0$.
\end{itemize}
Now, if we interlace $\sigma$ and $\rho$, then the sequence $\{x^k\}$ does not converge, it alternates between $x^0$ and $x^1$ indefinitely:
$
x^0 = 2 
\overset{\rho}{\longrightarrow}
x^1 = 2+12e^{-4} 
\overset{\sigma}{\longrightarrow}
x^0 = 2
\overset{\rho}{\longrightarrow}
x^1 = 2+12e^{-4} 
.
$
Thus, the convergence analysis of \MGOPT is incomplete, and as a side product, the method we propose establishes the convergence of \MGOPT as a special case 
(if gradient descent is used as the update step and our other assumptions hold), to be discussed in the contribution section.
We remark that the issue of non-convergence of \MGOPT has also been addressed using both line-search \cite{wen2010line} and trust region methods \cite{gratton2008recursive}.

\end{remark}

%
%
\subsection{Proximal gradient method}\label{sec:intro:subsec:prox}
Nowadays subgradient \cite{rockafellar1970convex} and proximal operator \cite{moreau1962fonctions} are standard tools for designing first-order algorithms to solve nonsmooth optimization problems \cite{combettes2005signal,parikh2014proximal,beck2017first},
especially for large-scale optimization where computing higher-order derivatives (e.g. the Hessian) is not feasible.
Here we give a quick review of the proximal operator and the proximal gradient operator.
We review subgradients in \cref{sec:mgprox:subsec:subdifferential}.

Rooted in the concept of Moreau's envelope \cite{moreau1962fonctions}, the proximal gradient method was first introduced in the 1980s in \cite[Eq. (4)]{fukushima1981generalized} as a generalization of the proximal point method \cite{rockafellar1976monotone}.
Under the abstraction of monotone operators, the proximal gradient method is understood as a forward-backward algorithm \cite{passty1979ergodic}, and it was later popularized by \cite{combettes2005signal} as the proximal forward-backward splitting. 
Nowadays proximal gradient method is ubiquitous in machine learning \cite{parikh2014proximal}.

The proximal gradient method solves problems of the form \eqref{prob:minfg} as follows.
Starting from an initial guess $x^1$, the method updates the variable by a gradient descent step (with a stepsize $\alpha > 0$) followed by a proximal step associated with $g_0$:
\begin{subequations}\label{eq:5}
\begin{align}
x^{k+1} &\,= \prox_{\alpha g_0} 
\Big( 
x^k - \alpha \nabla f_0(x^k)
\Big), 
\label{updt:proxgrad}
\\
\prox_{\alpha g_0}(x) 
& \coloneqq 
\underset{u}{\argmin} ~ \Big\{\, \alpha g_0(u) + \frac{1}{2} \| u - x \|_2^2 \,\Big\}
= \underset{u}{\argmin} ~ \Big\{\, g_0(u) + \frac{1}{2\alpha} \| u - x \|_2^2 \,\Big\}.
\label{def:prox_opertor}
\end{align}
\end{subequations}
If $f_0$ is $L_0$-smooth \eqref{def:Lsmooth}, we can set stepsize $\alpha$ in \eqref{updt:proxgrad} as  $\alpha \in (0, \frac{2}{L_0} )$ because such a stepsize brings strict functional decrease, and convergence to critical points 
\cite{beck2017first}.
The proximal operator \eqref{def:prox_opertor} itself is also an optimization problem, and in practice many commonly used $g_0$ are ``proximable'' in that \eqref{def:prox_opertor} has an efficiently computable closed-form solution.
The proximal gradient method has many useful properties.
To keep the introduction short, we introduce these properties later when needed.

%
%
\subsection{Contributions}
In this work, our contributions are:
\begin{enumerate}
\item We propose \MGProx (multigrid proximal gradient method) to solve \eqref{prob:minfg}.
It generalizes \MGOPT on smooth problems to nonsmooth problems using proximal gradient as the smoothing method.
A key ingredient in \MGProx is a newly introduced adaptive restriction operator in the multigrid process to handle the Minkowski sum of subdifferentials.
The key idea is about collapsing a set-valued vector into a singleton vector to ease computation, more to be explained in \cref{sec:mgprox}.

\item We provide theoretical results for 2-level \texttt{MGProx}: we show that
\begin{itemize}
\item \MGProx exhibits a fixed-point property, see Theorem~\ref{thm:fpp};

\item the coarse correction update (in step (v) of Algorithm~\ref{aglo:2levelMG_exact}) is a descent direction for the fine variable of the original fine level problem, see the subdifferential obtuse angle condition in Theorem~\ref{thm:coarse_dir_descent} and Lemma~\ref{lemma:alphaexists} for the existence of a coarse correction stepsize that provides a descent condition; 

\item the sequence $\{ f(x^k_0) \}_{k \in \IN}$ at the finest level converges to the optimal value  with rate $
1/k$ and $\big( 1 - \frac{\mu_0}{L_0} \big)^k$
, see Theorems~\ref{thm:rate1k} and \ref{thm:mgprox_converge_proxPL};
this result also establishes the convergence of \MGOPT (for $\sigma$ being the gradient update) in the convex case (see remark~\ref{remark:NashWrong}).

\item if we combine \MGProx with Nesterov's acceleration, the sequence $\{ f(x^k_0) \}_{k \in \IN}$ at the finest level converges to the optimal value with rate $1/k^2$, see Theorem~\ref{thm:FMGProx_fastr_ate}.
\end{itemize}

\item On the elastic obstacle problem, we show that multigrid  accelerates the proximal gradient method; we show that \MGProx runs faster than other methods.
See \cref{sec:exp}.

\end{enumerate}


%
%
\subsection{Literature review}\label{sec:intro:subsec:review}
The idea of multigrid is natural when handling large-scale elliptic PDE problems.

\subsubsection{Early works}
Early multigrid methods for non-smooth problems like \eqref{prob:minfg} pertain to the case of constrained optimization problems where $g_0$ is an indicator function on the feasible set.
For example, \cite{brandt1983multigrid} and \cite{mandel1984multilevel} develop multigrid methods for a symmetric positive definite (SPD) quadratic optimization problem with a bound constraint, which is equivalent to a linear complementarity problem. 
This applies, for example, to linear Elastic Obstacle Problems where $g_0$ is a box indicator function that models non-penetration constraints.
In \cite{hackbusch1983multi} this is extended to more general constrained nonlinear variational problems with SPD Fr\'echet derivatives, and to their associated nonlinear variational inequalities.
Later \cite{graser2009truncated} developed a Newton-MG (see below) method for an SPD quadratic optimization problem with more general but separable nonsmooth $g_0$.
This is extended in \cite{graser2019truncated} to a nonlinear objective function with nonsmooth $g_0$.

\subsubsection{Two families of multigrid}
We emphasize that there are at least two different approaches to perform multigrid in optimization.
The 2-level MGOPT algorithm (Algorithm \ref{algo:mgopt}) is an example of a \emph{full approximation scheme} (FAS) multigrid method for nonlinear problems.
The FAS approach, which was first described in \cite{brandt1977multi}, adds a $\tau$-correction term to the coarse nonlinear problem to ensure that the multigrid cycle satisfies a fixed-point property.
And because of the $\tau$, FAS is also called tau-correction method.

There is an alternative multigrid approach for solving nonlinear problems, which is the so-called Newton-multigrid method  (Newton-MG, \cite[Ch. 6]{briggs2000multigrid}), where the fine-level problem is first linearized using Newton's method and the linear systems in each Newton iteration are solved approximately using a linear multigrid method.
In other words, Newton-MG applies multigrid on solving the linear system $\nabla^2 f(x) d = \nabla f(x)$.
Newton-MG includes the works 
\cite{kornhuber1994monotone, graser2009truncated, gratton2010numerical, kocvara2016first, graser2019truncated}.

In the context of optimization problems, nonlinear multigrid methods can be devised to either work directly on the optimization problem and coarse versions of the optimization problem (as \MGOPT does), or they can be designed to work on the fine-level optimality conditions and coarse versions of them.

\subsubsection{How our approach differs}
Our method is a FAS approach like \cite{brandt1983multigrid,hackbusch1983multi,mandel1984multilevel} but our approach applies to general $g_0$ functions that go beyond indicator functions and include nonsmooth regularizations. 
While \cite{brandt1983multigrid, mandel1984multilevel} deal with linear problems, our approach applies to general nonlinear $f_0$. 
In contrast to \cite{brandt1983multigrid, hackbusch1983multi}, we don't use injection for the restriction operation, which often leads to slow multigrid convergence, but instead we use an adaptive restriction and interpolation mechanism that precludes coarse-grid updates to active points.

Our adaptive restriction and interpolation mechanism is similar to the truncation process used in 
\cite{graser2009truncated,graser2019truncated}, but our approach uses a FAS framework while \cite{graser2009truncated,graser2019truncated} use Newton-MG, and, most important we provide a convergence proof with convergence rates $1/k$, $(1-\mu_0/L_0)^k$ and $1/k^2$, while \cite{graser2019truncated} has no result on convergence rate.
Furthermore, Newton-MG requires the computation of 2nd-order information (the Hessian), while MGProx is a 1st-order method.

To sum up, our approach \textbf{is a first-order method} that avoids computing second derivative,
and the method \textbf{is a FAS} that does not require solving the equations $\nabla^2 f(x) d = \nabla f(x)$.
While existing multigrid methods in optimization are problem specific, our approach is \textbf{general for a class of non-smooth functions}.

\subsubsection{Multigrid outside PDEs}
\paragraph{Multigrid in image processing}
Besides PDEs, multigrid was used in the 1990s in image processing for solving problems with a nondifferentiable total variation semi-norm in image recovery (e.g., \cite{vogel1996iterative,chantextordmasculine1998multigrid}).
Note that these works bypassed the non-smoothness by \textit{smoothing} the total variation term, making them technically only solving \eqref{prob:mgopt} but not \eqref{prob:minfg}.

\paragraph{Multigrid in machine learning}
In the 2010s multigrid started to appear in machine learning, e.g., $\ell_1$-regularized least squares \cite{treister2012multilevel} and 
Nonnegative Matrix Factorization \cite{gillis2012multilevel}.
We remark that these works are not true multigrid method as there is no $\tau$ in the schemes, nor is the information of the fine variable carried to the coarse variable when solving the problem.

\paragraph{Recent work}
Recently \cite{parpas2017multilevel} proposed a multilevel proximal gradient method with a FAS structure, however it bypassed the technically challenging part of nonsmoothness by using smoothing, making it similar to \cite{vogel1996iterative,chantextordmasculine1998multigrid} in that they are only solving \eqref{prob:mgopt} but not \eqref{prob:minfg}.

The table below summarizes the comparison.
In the table, ``1st-order'' means the method discussed in the paper is a 1st-order method,
\\
\begin{tabular}{l|llll}
Work & 1st-order & FAS & convergence theory & general nonsmooth $g_0$ 
\\ \hline\hline
\cite{kornhuber1994monotone} & no & no & yes but no rate & no, box constraints only
\\
\MGOPT \cite{nash2000multigrid} &yes and no &  yes &  no (Remark\ref{remark:NashWrong})  & no 
\\
\cite{graser2009truncated} & no  & no &  no & yes
\\
\cite{gratton2010numerical} & no & no & no & no, box constraints only
\\
\cite{kocvara2016first} & yes & yes & yes but no rate & no, box constraints only
\\
\cite{parpas2017multilevel} & yes & yes & yes but smoothing & yes
\\
\cite{graser2019truncated} & no & no &  yes but no rate   & yes 
\\ \hline
\textbf{This work} & yes  & yes &  yes with rates &  yes
\end{tabular}

\subsection{Organization}
In \cref{sec:mgprox} we present a 2-level \MGProx and discuss its theoretical properties.
Then we present an accelerated \MGProx in \cref{sec:NestMGProx} and a multi-level \MGProx in \cref{sec:multilevel}.
In \cref{sec:exp} we demonstrate the performance of \MGProx compared with other methods.
We conclude the paper in \cref{sec:conc}.

\section{A two-level multigrid proximal gradient method}\label{sec:mgprox}
In \cref{sec:mgprox:subsec:RP}-\ref{sec:mgprox:subsec:adaptiveR}, 
we review subgradients for nonsmooth functions,
discuss their interaction with restriction (the coarsening operator), 
introduce the notion of adaptive restriction, and define the $\tau$ vector that carries the cross-level information.
We introduce a 2-level \MGProx method in \cref{sec:mgprox:subsec:2lv} and we provide theoretical results about the algorithm: 
fixed-point property (Theorem~\ref{thm:fpp}), 
descent property (Theorem~\ref{thm:coarse_dir_descent}),
existence of coarse correction stepsize (Lemma~\ref{lemma:alphaexists}) 
and convergence rates (Theorems~\ref{thm:rate1k} and \ref{thm:mgprox_converge_proxPL}).
In \cref{sec:mgprox:subsec:tau} we discuss further details of $\tau$.

%
%
\subsection{Functions at different levels}\label{sec:mgprox:subsec:RP}
Following \cref{sec:intro:subsec:mgopt}, 
we use $f_\ell : \IR^{n_\ell} \rightarrow \IR$ to denote functions at different coarse levels.
In this section we will focus on $\ell \in \{0,1\}$ but we remark that all the notations and definitions are generalized to $\ell \in \{0,1,\dots,L\}$ in Section~\ref{sec:multilevel}.
We denote the restriction of the fine objective function $F_0$ in \eqref{prob:minfg} as
$ F_1 \coloneqq \cR(F_0)  = \cR(f_0 + g_0)$, where $\cR$ is defined below.

\begin{definition}[Restriction]\label{def:restriction}
At a level $\ell \in \IN$, given a function $f_\ell : \IR^{n_\ell} \rightarrow \IR$, the restriction $\cR$ of $f_\ell$, denoted as $f_{\ell+1} \coloneqq \cR(f_\ell)$, is defined as 
$f_{\ell+1}(x_{\ell+1}) \coloneqq f_{\ell}(R x_{\ell})$, 
where $R : \IR^{n_\ell} \rightarrow \IR^{n_{\ell+1}}$ is a restriction matrix .
We also define the associated prolongation matrix 
$P : \IR^{n_{\ell+1}} \rightarrow \IR^{n_\ell}$ 
 as $P = cR^\top$ where $c>0$ is a predefined constant.
\end{definition}

\paragraph{Adaptive restriction and non-adaptive restriction}
We recall that a contribution of this work is the introduction of the adaptive restriction, to be discussed in \cref{sec:mgprox:subsec:adaptiveR} (see Definition\,\ref{eg:adaptiveR_L1}).
To differentiate the classical non-adaptive restriction (and the associated prolongation) from the adaptive version, we denote the non-adaptive restriction by $\overline{\cR}, \overline{R},\overline{\cP},\overline{P}$, and denote the adaptive one by $\cR, R,\cP,P$.
We remark that Definition~\ref{def:restriction} can be used for both versions of restriction and prolongation.
We give an example of $\overline{R},\overline{P}$ in \cref{sec:exp}.

%
%
\subsection{Review of subdifferential of nonsmooth functions}\label{sec:mgprox:subsec:subdifferential}
The subdifferential \cite{rockafellar1970convex}
is a standard framework used in convex analysis to deal with nondifferentiable functions.
A convex function $g(x) : \IRn \rightarrow \overline{\IR} \coloneqq \IR \cup \{+\infty\}$ is called nonsmooth if it is not differentiable for some $x$ in $\IRn$.
A point $q \in \IRn$ is called a subgradient of $g$ at  $x$ if for all $y \in \IRn$ the inequality $g(y) \geq g(x) + \langle q, y -x \rangle$ holds.
The subdifferential of $g$ at a point $x$ is defined as the set of all subgradients of $g$ at $x$, i.e.,
\begin{equation}
  \partial g(x)
\coloneqq 
\Big\{~
q \in \IRn ~~\big|~~ 
g(y)  \geq g(x) + \langle q, y -x \rangle
~\forall y \in \IRn
~\Big\} ~ \subset \IRn,
\label{def:subdifferential_II}
\end{equation}
so $\partial g(x)$ is generally set-valued.
If $g$ is differentiable at $x$, then $\nabla g(x)$ exists and the set $\partial g(x)$ reduces to the singleton $\big\{\nabla g(x)\big\}$.

\paragraph{Subdifferential sum rule (Moreau–Rockafellar theorem)}~~
Let $\oplus$ denote the Minkowski sum.
Since subdifferentials are generally set-valued, hence for two functions $f, g$, generally $\partial (f+g) \neq \partial f \oplus \partial g$ but $\partial (f+g) \supset \partial f \oplus \partial g$.
The sum rule $\partial (f+g) = \partial f \oplus \partial g$ holds if $f,g$ satisfy a qualification condition (e.g. \cite[Theorem 3.36]{beck2017first}): the relative interior of the domain of $f$ has a non-empty intersection with the relative interior of domain of $g$,
i.e.,
\begin{equation}
\ri (\dom f) \,\cap\, \ri (\dom g) \neq \varnothing
\implies
\partial \Big( f(x) + g(x) \Big) =
\partial f(x) \oplus \partial g(x),
~~~
\forall x \in \textrm{dom} f \,\cap\, \textrm{dom} g.
\label{fact:subgrad_sum_rule}
\end{equation}
The right-hand side (RHS) of \eqref{fact:subgrad_sum_rule}, which is the subdifferential sum rule,
is known as the Moreau–Rockafellar theorem \cite{lai1988moreau}.
We now discuss the fact that the functions $f_\ell, g_\ell$ for all levels $\ell$ in this work satisfy the Moreau–Rockafellar theorem. 
\begin{itemize}
\item At level $\ell=0$, we have the left-hand side (LHS) of \eqref{fact:subgrad_sum_rule} for $f_0, g_0$ 
by assumption.
\item At levels $\ell>0$, by Definition~\ref{def:restriction}, 
the LHS of \eqref{fact:subgrad_sum_rule} holds for the coarse functions.
To see this, take the coarse function as composition of the fine function with a linear map, and recognize the fact that the domain of a function is preserved under composition with a linear map.
To be explicit, we have
\[
\begin{array}{rcl}
\ri \Big( \dom f_1 \Big)
\cap 
\ri \Big( \dom g_1 \Big)
~=~
\ri \Big( \underbrace{\dom (f_0 \circ R)}_{\text{whole }\IR^{n_1}} \Big)
\cap 
\ri \Big( \dom (g_0 \circ R)  \Big)
~\neq~ \varnothing.
\end{array}
\]
\end{itemize}
To sum up, in this work \eqref{fact:subgrad_sum_rule} holds for all levels $\ell$:
\begin{equation}
\partial F_\ell(x_\ell) 
\coloneqq  \partial \Big( f_\ell(x_\ell) + g_\ell(x_\ell) \Big)
= \partial f_\ell(x_\ell) \oplus \partial g_\ell(x_\ell)
= \nabla f_\ell(x_\ell) + \partial g_\ell (x_\ell)
,
\label{ass:F_sum_rule}
\end{equation}
where $+$ is used instead of $\oplus$ in $\nabla f_\ell(x_\ell) + \partial g_\ell (x_\ell)$ because $\nabla f_\ell(x_\ell)$ is a singleton.

\subsection{Convexity and subdifferential of coarse function}
From Definition~\ref{def:restriction}, the coarse function can be written as $F_{\ell+1} \coloneqq F_\ell \circ R$, meaning that we can see the restriction process as the fine function taking composition with the linear map $R$.
Such composition view point gives us a series of useful properties for this work.
First, we have a closed-form expression for the subdifferential of the coarse function in terms of the the subdifferential of the fine function.
I.e., by \cite[Theorem 3.43]{beck2017first}, we have that
\begin{equation}\label{subdiff_coarse_linearmap}
\partial F_{\ell+1}(x) = \partial (F_\ell \circ R)(x) = R^\top \partial F_\ell (Rx). 
\end{equation}
Then, by the fact that convexity is preserved under linear map, we have that ``restriction preserves convexity''.
In other words, $f_{\ell+1}$ is convex if $f_{\ell}$ is convex.
Furthermore, if $R$ is full-rank (which is the case in this paper), restriction is submultiplicative on the modulus of convexity, as illustrated in the following lemma.

\begin{lemma}[Composition with full-rank matrix preserves  convexity]
\label{lem:strcvx_linearmap}
Given a function $F : \IRn \rightarrow \overline{\IR}$ that is $\mu$-strongly convex and a rank-$n$ matrix $R \in \IR^{m \times n}$, the function $F \circ R$ is $\mu \sigma^2_n$-strongly convex, where $\sigma_n$ is the $n$th singular value of $R$.
\begin{proof}
$F$ is strongly convex so $\partial F$ is strongly monotone \cite{bauschke2011convex}: for all $x,y \in \dom F$, 
\begin{equation}\label{lem:strcvx_linearmap:fineF}
\big\langle x - y, ~ \partial F(x) - \partial F(y ) \big\rangle 
~\geq~
\mu \| x - y\|_2^2.
\end{equation}
Now we show $\partial (F \circ R)$ is also strongly monotone.
For all $x,y $ in $\dom (F \circ R)$, we have
\[
\begin{array}{rcl}
\langle 
x - y, 
~ \partial (F \circ R)(x) - \partial (F \circ R)(y )
\rangle 
&\overset{\eqref{subdiff_coarse_linearmap}}{=}&
\langle 
x - y, 
~ R^\top  \partial F(Rx) - R^\top  \partial F (Ry )
\rangle 
\\
&=&
\langle 
Rx - Ry, 
~ \partial F(Rx) -  \partial F (Ry )
\rangle 
\\
&\overset{\eqref{lem:strcvx_linearmap:fineF}}{\geq}&
\mu \| Rx-Ry\|_2^2
~\geq~ \mu \sigma_n^2 \| x-y\|_2^2,
\end{array}
\]
where $\sigma_n > 0 $ is the $n$th singular value of the full-rank matrix $R$.
Thus the subdifferential $\partial (F \circ R)$ is  is $\mu \sigma_n^2$-strongly monotone, and thus the function $F \circ R$ is $\mu \sigma_n^2$-strongly convex.
\end{proof}
\end{lemma}

\paragraph{Notation for sets}
From now on, when we encounter an expression containing both set-valued vector(s) and singleton vector(s), we underline the set-valued term(s) for visual clarity.

%
%
\subsection{Adaptive restriction and the \texorpdfstring{$\tau$}{} vector}\label{sec:mgprox:subsec:adaptiveR}
Since subdifferentials are set-valued, we define $\tau$ in \MGProx as an element of a set.
At a level $\ell$, we define a set $\underline{\tau_{\ell \rightarrow \ell+1}} ~\coloneqq~
\underline{\partial F_{\ell+1}(x_{\ell+1})} \oplus (- R) \underline{\partial F_\ell (x_\ell)} $, where $\ell \rightarrow \ell+1$ specifies that $\tau$ connects level $\ell$ to level $\ell+1$, and the matrix $R$ here is an adaptive restriction operator that we will define soon.
In \MGProx we choose an element of $\underline{\tau_{\ell \rightarrow \ell+1}}$ as the tau vector.
That is, at level $\ell=0$,
\begin{subequations}
\begin{align}
\tau_{0 \rightarrow 1}
~\, \in \,~ 
\underline{\tau_{0 \rightarrow 1}}
& 
~\coloneqq~
\underline{\partial F_1(x_1)} \oplus (- R) \underline{\partial F_0(x_0)}
\label{def:tau_mgprox1}
\\
&\overset{\eqref{ass:F_sum_rule}}{=} \nabla f_1(x_1) - R \nabla f_0(x_0) + \underline{\partial g_1(x_1)}
\oplus 
(- R) \underline{\partial g_0(x_0)}.
\label{def:tau_mgprox2}
\end{align}
\end{subequations}
Note that $\underline{\tau_{0 \rightarrow 1}}$ is a function of two points at two different levels.
In \eqref{def:tau_mgprox1} $\underline{\tau_{0 \rightarrow 1}}$ is the Minkowski sum of two subdifferentials $\partial F_1(x_1) $ and $-R\partial F_0(x_0)$ which are generally set-valued.
To obtain a tractable coarse-grid optimization problem (corresponding to line (iv) in Algorithm~\ref{algo:mgopt}) we need to avoid complications coming from the Minkowski sum, and we do this by modifying the standard restriction (and prolongation) 
by zeroing out columns in $\overline{R}$ to form $R$ such that the second subdifferential $R \partial g_0(x_0)$ in \eqref{def:tau_mgprox2}
is a singleton vector.
Similarly, we  zero out the corresponding rows in $\overline{P}$ to form $P$ for the coarse correction step, such that non-differentiable fine points are not corrected by the coarse grid.
This zeroing out process is adapted to the current point $x_0$, so we call this $R$ an adaptive restriction operator. 
In other words, the purpose of the adaptive restriction is to reduce a generally set-valued subdifferential $R\partial g_0(x_0)$ to a singleton.
We denote the adaptive operator $R$ corresponding to a point $x$ as $R(x)$ and thereby the adaptive restriction of $x$ is denoted as $R(x)x$.
Sometimes we just write $Rx$ if the meaning is clear from the context.
Based on the above discussion, we now formally define adaptive restriction,
and we give an example in \cref{sec:exp}.

\begin{definition}[Adaptive restriction operator for separable $g$]\label{eg:adaptiveR_L1}
For a possibly nonsmooth function $g : \IRn \rightarrow \IR$ that is separable, i.e., with $x = [x_1,x_2,\dots,x_n]$, $g(x) = \sum g_i(x_i)$ where $g_i$ is a function only of component $x_i$, 
given a full restriction operator $\overline{R}$ and a vector $x$, the adaptive restriction operator $R$ with respect to a function $g$ at $x$ is defined by zeroing out the columns of $\overline{R}$ corresponding to the elements in $\partial g$ that are set-valued.
\label{def:adaptiveR}
\end{definition}

\paragraph{The coarse problem is nonsmooth and $\tau$ is an element of a set}~~
Now it is clear that the subdifferential $R \partial g_0(x_0)$ in \eqref{def:tau_mgprox2} is a singleton.
From the fact that the coarse problem is nonsmooth (where the function $g_1$ is nonsmooth), there are two consequences:
\begin{itemize}
    \item it makes the coarse problems difficult to solve as well as the original problem.
    This makes our approach differ from works such as \cite{parpas2017multilevel} where the coarse problems are replaced by a smooth approximation;
    and

    \item the first subdifferential $\partial g_1(x_1)$ in \eqref{def:tau_mgprox2} is possibly set-valued, thus the RHS of \eqref{def:tau_mgprox2} is generally set-valued and so is $\underline{\tau_{0 \rightarrow 1}}$, and we define $\tau_{0 \rightarrow 1}$ to be a member of the set $\underline{\tau_{0 \rightarrow 1}}$.
We emphasize that in the algorithm to be discussed below we can pick any value for  $\tau_{0 \rightarrow 1}$ in the set.
We will explore the choice of $\tau$ after we have given a complete picture of \texttt{MGProx}.
\end{itemize}

\paragraph{Adaptive restriction differs from Kornhuber's basis truncation}
In the PDE literature there is a multigrid method called Kornhuber's basis truncation \cite{kornhuber1994monotone}, in which at first glance looks similar to \texttt{MGProx}.
We remark that Kornhuber's basis truncation is designed only for box-constrained optimization.
The truncation zeros out the basis of the optimization variable, while adaptive restriction zeros the subdifferential vector (see \eqref{def:tau_mgprox2}).
Also, the truncation is only applied to the finest level \cite[Section 2.2]{kocvara2016first}, while \MGProx applies adaptive restriction applies to all the levels.

We are now ready to present \texttt{MGProx}.
Here we present a 2-level \MGProx method for illustration, and we move to a general multi-level version in \cref{sec:multilevel}.
For adaptive $R$, now all the Minkowski additions are trivial addition so we use $+$ instead of $\oplus$.


%
%
\subsection{A 2-level \texorpdfstring{\MGProx}{}algorithm}\label{sec:mgprox:subsec:2lv}
Similar to the 2-level \MGOPT method for solving Problem \eqref{prob:mgopt}, we propose a 2-level \MGProx method (Algorithm \ref{aglo:2levelMG_exact}) that solves Problem \eqref{prob:minfg} by utilizing a coarse problem defined as 
\begin{equation}\label{prob:minfg_coarse}
\argmin_{\xi \in \IR^{n_1}}\,
\bigg\{\,
F^{\tau}_1 (\xi) 
~\coloneqq~  F_1(\xi) - \langle \tau_{0 \rightarrow 1}^{k+1}, \xi \rangle
~=~ f_1 (\xi) + g_1(\xi) - \langle \tau_{0 \rightarrow 1}^{k+1}, \xi \rangle
\,\bigg\}.
\end{equation}

\begin{algorithm}
\caption{2-level \MGProx for an approximate solution of \eqref{prob:minfg}}
\label{aglo:2levelMG_exact}
\begin{algorithmic}
\STATE{Initialize $x^1_0$, $R$ and $P$}
\FOR{$k = 1,2,\dots$}
\STATE{(i)~~~~ $y^{k+1}_0 = \prox_{\frac{1}{L_0}g_0} \Big (x^k_0 - \frac{1}{L_0} \nabla f(x^k_0) \Big)$
\hfill level-0 proximal gradient step \vspace{0.25mm}
}
\STATE{(ii)~~~ $y^{k+1}_1  = R(y^{k+1}_0) y^{k+1}_0 $   
\hfill construct the level-1 coarse variable \vspace{0.25mm}
}
\STATE{(iii)~~ $\tau_{0 \rightarrow 1}^{k+1} \hspace{-1mm} \in\,
\underline{\partial F_1(y^{k+1}_1 )}  - R(y^{k+1}_0) \, \underline{\partial  F_0(y^{k+1}_0)}$
\hfill construct the tau vector \vspace{0.25mm}
}
\STATE{(iv)~~ $
x^{k+1}_1
= \underset{\xi}{\argmin} \, \Big\{ F_1^{\tau}(\xi) \coloneqq F_1(\xi) - \langle \tau_{0 \rightarrow 1}^{k+1}, \xi \rangle \Big\}
$
\hfill solve the level-1 coarse problem
}
\STATE{(v)\,~~~ $z^{k+1}_0 = y^{k+1}_0  + \alpha P\big(  x^{k+1}_1 - y^{k+1}_1 \big)$
\hfill coarse correction \vspace{0.25mm}
}
\STATE{(vi)~~\, $x^{k+1}_0 = \prox_{\frac{1}{L_0}g_0} \Big (z^{k+1}_0 - \frac{1}{L_0} \nabla f(z^{k+1}_0) \Big)$
\hfill level-0 proximal gradient step~~
}
\ENDFOR
\end{algorithmic}
\end{algorithm}

Here are remarks for the steps in Algorithm~\ref{aglo:2levelMG_exact}.
\begin{itemize}
    \item (i): we perform one or more proximal gradient iterations on the fine variable with a constant stepsize $1/L_0$, where $L_0$ is the Lipschitz constant of $\nabla f_0$.
    \item (iii): we pick a value within the set to define $\tau$;
    as we are now using adaptive $R$, we use $+$ instead of $\oplus$ in the expression of $\tau$.
    \item (v): $\alpha > 0$ is a stepsize; for its selection see \cref{sec:mgprox:subsec:alpha}.
    \item The restriction for variable $y_0^{k+1}$ and the restriction for the subdifferential $ \underline{\partial  F_0(y^{k+1}_0)}$ can be slightly different.
    On $y_0^{k+1}$ the restriction is the full restriction, on $ \underline{\partial  F_0(y^{k+1}_0)}$ is the adaptive one. 
    The explanation is as follows.   
    For the particular cases of $g_0$ such as $\ell_1$, $\max\{ \cdot, 0\}$ (element-wise maximum) and $\iota_{[0,\infty)}$ (indicator of nonnegative orthant), when we zero out column $i$ of $R$, the corresponding entry of $y_0^{k+1}$ is already $0$.    
    Note that this conclusion does not hold in general for other nonsmooth functions or when the non-differentiability occur at another point (say at $x=1$). 
    Generally in those cases we will need to specify which restriction matrix (the full $R$ or the adaptive $R$) to be used to define $P$.
        However, to simplify the presentation we always assume the non-differentiability occurs at $x=0$,
         by performing a translation to shift the non-differentiability to occur at $x=0$.
\end{itemize}

%
%
\subsubsection{Fixed-point property}
Algorithm~\ref{aglo:2levelMG_exact} exhibits the following fixed-point property.
%
%
\begin{theorem}[Fixed-point]\label{thm:fpp}
In Algorithm~\ref{aglo:2levelMG_exact}, if $x^k_0$ solves \eqref{prob:minfg}, then we have the fixed-point properties $x^{k+1}_0 =y^{k+1}_0 = x^k_0$ and $x^{k+1}_1 = y^{k+1}_1$.
\begin{proof}
The fixed-point property of the proximal gradient operator \cite[page 150]{parikh2014proximal} gives
\begin{equation}\label{thm:fpp:e1}
    y^{k+1}_0 \overset{\textrm{fixed-point}}{=}  x^k_0  
    \overset{\textrm{assumption}}{=}  \argmin \, F_0(x).
\end{equation}
As a result, the coarse variable satisfies
$y^{k+1}_1 \coloneqq R y^{k+1}_0  \overset{\eqref{thm:fpp:e1}}{=} Rx^k_0
$.
The 1st-order optimality of 
$ y^{k+1}_0 
\overset{\eqref{thm:fpp:e1}}{=} \argmin F_0
$ gives $0 \in 
\underline{\partial F_0( y^{k+1}_0 )}$.
Multiplying by $-R$ (which reduces the set $\underline{\partial F_0(x^k_0)}$ to a singleton) gives 
\begin{equation}\label{thm:fpp:e3}
 0 = - R \underline{\partial F_0(x^k_0)}.
\end{equation}
Then adding $\partial F_1( y^{k+1}_1)$ on both sides of \eqref{thm:fpp:e3} gives
\begin{subequations}
\begin{align}
&&\underline{\partial F_1( y^{k+1}_1)}  &= \underline{\partial F_1( y^{k+1}_1)} - R(x^k_0) \underline{\partial F_0(x^k_0) }
\label{thm:fpp:e4a}
\\
&&                          &\overset{\eqref{def:tau_mgprox1}}{\ni} \tau_{0 \rightarrow 1}^{k+1} 
\label{thm:fpp:e4b}
\end{align}
\end{subequations}
In \eqref{thm:fpp:e3}, $- R \underline{\partial F_0(x^k_0)} $ is the zero vector, so the equality in \eqref{thm:fpp:e4a} holds since we are adding zero to a (non-empty) set.
The inclusion \eqref{thm:fpp:e4b} follows from \eqref{def:tau_mgprox1} as the expression $\underline{\partial F_1( y^{k+1}_1)} - R(x^k_0) \underline{\partial F_0(x^k_0) }$ is the set $\underline{\tau_{0 \rightarrow 1}^{k+1}}$ by definition.

Now rearranging \eqref{thm:fpp:e4b} gives $0 \in \underline{\partial F_1( y^{k+1}_1)}  - \tau_{0 \rightarrow 1}^{k+1}$,
which is exactly the 1st-order optimality condition for the coarse problem
$ \underset{\xi}{\argmin} \, F_1( \xi ) -  \big\langle \tau_{0 \rightarrow 1}^{k+1} , \xi \big\rangle$.
By strong convexity of $ F_1( \xi ) -  \big\langle \tau_{0 \rightarrow 1}^{k+1} , \xi \big\rangle$, the point $y^{k+1}_1$ is the unique minimizer of the coarse problem, so $x^{k+1}_1 = y^{k+1}_1$ by step (iv) of the algorithm and $x^{k+1}_0 = y^{k+1}_0 \overset{\eqref{thm:fpp:e1}}{=} x^k_0$ by steps (v) and (vi).
\end{proof}
\end{theorem}

Theorem~\ref{thm:fpp} shows that at convergence,
we have fixed-point $x_0^{k+1}=y_0^{k+1}$ at fine level and also $x_1^{k+1} = y_1^{k+1}$ at the coarse level.
Next we show that when $x_1^{k+1} \neq y_1^{k+1}$,
the objective function value sequence is converging.

%
%
\subsubsection{Coarse correction descent: angle condition}
In nonsmooth optimization, descent direction properties are drastically different from smooth optimization \cite{noll2014convergence}.
For example for the subgradient method, the classical angle condition no longer describes a useful set of search directions for the subgradient.
In \MGProx the coarse correction direction $P(x^{k+1}_1 - y^{k+1}_1)$ is a nonsmooth descent direction, and we will show that $P(x^{k+1}_1 - y^{k+1}_1)$ decreases the objective function value, based on the theorem below and Lemma~\ref{lemma:alphaexists}.
%
%
\begin{theorem}[Angle condition of coarse correction]\label{thm:coarse_dir_descent}
If $P(x^{k+1}_1 - y^{k+1}_1) \neq 0$, then
\begin{equation}\label{iq:descent_dir}
\big\langle 
\underline{\partial F_0(y^{k+1}_0)} , P(x^{k+1}_1 - y^{k+1}_1) 
\big\rangle 
<
0. 
\end{equation}
\end{theorem}
Before we prove the theorem we emphasize that \eqref{iq:descent_dir} applies for any subgradient in the set $\underline{\partial F_0(y^{k+1}_0)} $.
Furthermore,
\[
\eqref{iq:descent_dir} 
\iff
\big\langle 
P^\top \underline{\partial F_0(y^{k+1}_0)} , x^{k+1}_1 - y^{k+1}_1 
\big\rangle
< 0
\overset{P^\top = cR}{\iff}
c \big\langle 
R \underline{\partial F_0(y^{k+1}_0)} , x^{k+1}_1 - y^{k+1}_1 
\big\rangle
< 0
.
\]
As $c>0$, showing \eqref{iq:descent_dir} is equivalent to showing
\begin{equation}
\big\langle R \underline{\partial F_0(y^{k+1}_0)},  ~ x^{k+1}_1 - y^{k+1}_1 \big\rangle < 0,
\label{pf:descent_to_prove}
\end{equation}
where $R \underline{\partial F_0(y^{k+1}_0)} $ is a singleton vector for all subgradients in $\underline{\partial F_0(y^{k+1}_0)}$ due to the adaptive $R$.
\begin{proof}
By definition 
$\tau_{0 \rightarrow 1}^{k+1} \overset{\eqref{def:tau_mgprox1}}{\in} \underline{\partial F_1(y_1^{k+1})}  - R \underline{\partial F_0(y_0^{k+1})}$ 
and the fact that $R \underline{\partial F_0(y^{k+1}_0)} $ is a singleton, we have
$
R \underline{\partial F_0(y_0^{k+1})} 
\in
\underline{\partial F_1(y_1^{k+1})}  - \tau_{0 \rightarrow 1}^{k+1}
\overset{\eqref{prob:minfg_coarse}}{=} \underline{\partial F_1^\tau(y_1^{k+1})}$,
showing that $R \underline{\partial F_0(y_0^{k+1})} $ is a subgradient of $F_1^{\tau}$ at $y_1^{k+1}$.
For any subgradient in the subdifferential $\underline{\partial F_1^\tau(y_1^{k+1})} $, we have the following which implies \eqref{pf:descent_to_prove}:
$
\big\langle \underline{\partial F_1^{\tau}(y^{k+1}_1)}, x^{k+1}_1 - y^{k+1}_1 \big\rangle
< F_1^{\tau}(x^{k+1}_1) - F_1^{\tau}(y^{k+1}_1) 
< 0,
$
where the first strict inequality is due to $F_1^\tau$ being a strongly convex function (which implies strict convexity);
the second inequality is by $x^{k+1}_1 \coloneqq \underset{\xi}{\argmin}\,
F_1^{\tau}(\xi)$ and the assumption that $x_1^{k+1} \neq y_1^{k+1}$.
\end{proof}

\begin{remark}
Theorem~\ref{thm:coarse_dir_descent} holds 
for convex (not strongly convex)  $f_0$
by changing $<$ to $\leq$.
\end{remark}

%
%
\subsubsection{Existence of coarse correction stepsize \texorpdfstring{$\alpha_k$}{}}
\label{sec:mgprox:subsec:alpha}
Based on Theorem~\ref{thm:coarse_dir_descent}, we now show that there exists a stepsize $\alpha_k > 0$ such that 
\begin{equation}
F_0(z^{k+1}_0) 
~\coloneqq~
F_0\big(
y^{k+1}_0 + \alpha_k P(x^{k+1}_1 - y^{k+1}_1) 
\big)
~<~
F_0(y^{k+1}_0).
\label{iq:coarse_descent_on_F}
\end{equation}

\begin{lemma}[Existence of stepsize]\label{lemma:alphaexists}
There exists $\alpha_k > 0$ such that \eqref{iq:coarse_descent_on_F} is satisfied for $P(x^{k+1}_1 - y^{k+1}_1) \neq 0$.
\end{lemma}

To prove the lemma, we make use 
of a fact of the subdifferentials of finite convex functions \cite[Def.1.1.4, p.165]{hiriart2004fundamentals}: $\partial g(x)$ is a nonempty compact convex set $\cS \subset \IRn$ whose support function $\displaystyle \sup_{s}\{ \langle s, x\rangle \,|\, s \in \cS  \}$ 
is the directional derivative of $g$ at $x$.
By this, the subdifferential
 $\underline{\partial F_0(y_0^{k+1})}$ is a compact convex set whose support function is the directional derivative of $F_0$ at $y_0^{k+1}$.
We emphasize that here 
$F_0(y_0^{k+1})$ is 
finite, 
so we can make use of the result on directional derivative in \cite[Def. 1.1.4, p.165]{hiriart2004fundamentals}, which only applies for finite convex functions, associated with the subdifferential.

\begin{proof}
We prove the lemma in 3 steps.
\begin{enumerate}

\item (Half-space)
The strict inequality in Theorem~\ref{thm:coarse_dir_descent} means that $\underline{\partial F_0(y_0^{k+1})}$ is strictly inside a half-space with normal vector $p = P(x^{k+1}_1 - y^{k+1}_1)$.  

\item (Strict separation)
Being a compact convex set, $\underline{\partial F_0(y_0^{k+1}0)}$ lying strictly on one side of the hyperplane must be a positive distance (say $\alpha_k > 0$) from that hyperplane.

\item (Support and directional derivative)
Evaluating the support function of $\underline{\partial F_0(y_0^{k+1})}$, i.e., the directional derivative of $F_0$ at $y_0^{k+1}$ in the direction $p$, we have \eqref{iq:coarse_descent_on_F}.
\end{enumerate}
\end{proof}

\begin{remark}[On the compactness of subdifferential]
\label{remark:compact}
For a function $\phi$, the set $\partial \phi$ is compact on $\textrm{int}~\dom \phi$.
Note that $\partial \phi$ is not compact for indicator functions at the boundary. 
So, for Lemma~\ref{lemma:alphaexists} to hold for function $g_0$, we assume $\dom g_0 \coloneqq \IRn$ with $g_0 : \IRn \rightarrow \IR$ as in the \textbf{Assumption} in the introduction.  
Such an assumption is needed for the proof that a positive $\alpha$ exists satisfying the line-search condition.
The impact of this assumption is that we are not allowing $g_0$ to be an indicator function, and thus reducing the scope of the applicability of the theory of \texttt{MGProx}.
However,
\begin{itemize}
    \item Empirically, we have observed that \MGProx works for $g_0$ being an indicator function (such as box constraints).
    
    \item The convergence proof of \MGProx does not require a positive coarse correction stepsize $\alpha$ (i.e., allowing $\alpha < 0$), 
    so in principle we can relax the condition of $g_0 : \IRn \rightarrow \IR$ to $g_0 : \IRn \rightarrow \overline{\IR}$.
    
    \item It will be an interesting future work 
    to generalize Lemma~\ref{lemma:alphaexists} for $g_0 : \IRn \rightarrow \overline{\IR}$.
    There is a different proof showing that a positive $\alpha$ exists in the case of certain indicator functions. 
    For example, the following is a theorem.
    If $H$ is a hyperplane and $P$ is a polyhedral set (closed but possibly non-compact) lying strictly on one side of $H$, then there is a positive distance between $H$ and $P$.
\end{itemize}
 
\end{remark}

Now we see that Theorem~\ref{thm:coarse_dir_descent} implies Lemma~\ref{lemma:alphaexists} which then implies the descent condition \eqref{iq:coarse_descent_on_F}.
Later in Theorem~\ref{thm:rate1k} and Theorem~\ref{thm:mgprox_converge_proxPL}, by using \eqref{iq:coarse_descent_on_F} together with the sufficient descent property of proximal gradient (\cref{lemma:suff_descent_proxgrad} below), we prove that the sequence $\big\{F_0(x^k_0)\big\}_{k \in \IN}$ produced by Algorithm~\ref{aglo:2levelMG_exact} converges to the optimal value $F^* \coloneqq \min F$.
Before that, in the next paragraph we first discuss about tuning the coarse correction stepsize.

\begin{lemma}
(Sufficient descent property of proximal gradient \cite[Lemma 10.4]{beck2017first})
\label{lemma:suff_descent_proxgrad} 
Let $L_0$ be the Lipschitz constant of $\nabla f_0$.
For step (i) in Algorithm~\ref{aglo:2levelMG_exact},  we have
\begin{equation}
F_0( y^{k+1}_0)
\leq 
F_0( x^k_0) - \| G_0(x_0^k) \|_2^2 / (2L_0),
~~
G_0( x_0^k) 
=
L_0 \big[
x_0^k - \prox_{ \frac{1}{L} g_0 }\big( x_0^k - \nabla f_0(x_0^k)/L_0
\big)
\big].
\label{eqn:suffdescentPGD}
\end{equation}
\end{lemma}
$G_0( x_0^k)$ is called the proximal gradient map of $F_0$ at $x_0^k$.
The inequality also holds for step (vi).

\begin{remark}[Theorem~\ref{thm:coarse_dir_descent} holds for inexact coarse update]
Step (iv) in Algorithm~\ref{aglo:2levelMG_exact} can be expensive.
If we replace (iv) by an iteration of (coarse level) proximal gradient step, then by \eqref{eqn:suffdescentPGD} we have $F_1^\tau(x_1^{k+1}) \leq F_1^\tau(y_1^{k+1}) - \| G_1^\tau(y_1^{k+1})\|_2^2 / (2L_1)$.
For $x_1^{k+1} \neq y_1^{k+1}$ we have $F_1^\tau(x_1^{k+1}) - F_1^\tau(y_1^{k+1}) < 0$.
Thus the descent condition holds for an inexact coarse update.
\end{remark}

%
%
\subsubsection{Tuning the coarse correction stepsize \texorpdfstring{$\alpha_k$}{}}
\label{sec:mgprox:subsec:alpha_tune}
First, exact line search is impractical: finding $\alpha_k \coloneqq 
\underset{\alpha \geq 0}{\argmin}
F_0\big( 
y^{k+1}_0 + \alpha P (x^{k+1}_1 - y^{k+1}_1) 
\big)
$
is generally expensive.
Next, classical inexact line searches such as the Wolfe conditions, Armijo rule, Goldstein line search (e.g., see \cite[Chapter 3]{nocedal1999numerical}) cannot be used here as they were developed for smooth functions.
While it is possible to develop nonsmooth version of these methods, such as a nonsmooth Armijio rule in tandem with backtracking on functions that satisfy the Kurdyka-\L ojasiewicz inequality with other additional conditions in \cite{noll2014convergence}, this is out of the scope of this work. 
Precisely, consider the condition
$
F_0\big( 
y^{k+1}_0 + \alpha P (x^{k+1}_1 - y^{k+1}_1) 
\big) 
\leq 
F_0\big( y^{k+1}_0 \big)
+ c_1 \langle \underline{\partial F_0(y_0^{k+1})}, \alpha P (x^{k+1}_1 - y^{k+1}_1) \rangle
$
where $c_1 \in (0,1)$.
Due to the strict inequality \eqref{iq:descent_dir}, the value $c_1$ is possibly fluctuating and there is no  simple-and-efficient way to determine its value. 
Thus, for this paper, we use simple naive backtracking as shown in Algorithm~\ref{aglo:naivelinesearch}, which just enforces \eqref{iq:coarse_descent_on_F} without any sufficient descent condition.
While we acknowledge that the traditional wisdom in optimization tells that naive descent conditions such as \eqref{iq:coarse_descent_on_F} are generally not enough to obtain convergence to the optimal point, we note that \MGProx is not solely using the coarse correction to update the variable; instead it is a chain of interlaced iterations of proximal gradient descent and coarse correction, and we will show next that the sufficient descent property of proximal gradient descent \eqref{eqn:suffdescentPGD} alone provides enough descending power for the function value $F_0$ to convergence to the optimal value.

\begin{algorithm}
\caption{Naive line search}\label{aglo:naivelinesearch}
\begin{algorithmic}
\STATE{Set $\alpha > 0$ and select a tolerance $\epsilon > 0$  (e.g. $10^{-15}$)}
\WHILE{true}
\STATE{\textbf{If} $F_0\big( y^{k+1}_0 + \alpha P (x^{k+1}_1 - y^{k+1}_1) \big) \leq F_0( y^{k+1}_0 )$ \textbf{ then }
return $z_0^{k+1}=y^{k+1}_0 + \alpha P(x^{k+1}_1 - y^{k+1}_1)$.}
\\
\STATE{\textbf{else if} $\alpha > \epsilon $ \textbf{then} $\alpha = \alpha/2 $.}
\\
\STATE{~~~~~~~\textbf{else} return $z_0^{k+1}=y^{k+1}_0$.}
\ENDWHILE
\end{algorithmic}
\end{algorithm}


\begin{remark}[Algorithm~\ref{aglo:naivelinesearch} has a fixed complexity]
Suppose the algorithm runs $N$ iterations, then with a tolerance $\epsilon = 10^{-15}$, the coarse stepsize $\alpha = 1/2^N > \epsilon $ gives $N = 50 $.
I.e., Algorithm~\ref{aglo:naivelinesearch} will run at most 50 iterations to search for $\alpha$.
\end{remark}

\begin{remark}[Infinitesimal $\alpha$ may not be computable]
Lemma \ref{lemma:alphaexists} on the existence of $\alpha > 0$ does not exclude the possibility that $\alpha$ is too close to zero so that $z^{k+1}_0$ cannot be distinguished from $y^{k+1}_0$.
In this case Algorithm~\ref{aglo:naivelinesearch} simply gives $\alpha = 0$.
\end{remark}

%
%
\subsubsection{\texorpdfstring{$\cO(1/k)$}{} convergence rate}
Inequality \eqref{iq:coarse_descent_on_F} implies that in the worst case the coarse correction $P(x^{k+1}_1 - y_1^{k+1})$ in the multigrid process is ``doing nothing'' on $y_0^{k+1}$, which occurs when $P(x^{k+1}_1 - y_1^{k+1}) = 0$ or $x^{k+1}_1 = y^{k+1}_1$.
We now show that the descent inequality $F_0(x^{k+1}_0) \leq F_0(z^{k+1}_0) \leq F_0(y^{k+1}_0)$ implies that the convergence rate of the sequence $\big\{ F_0(x_0^k) \big\}_{k \in \IN} $ for $\{x_0^{k} \}_{k\in \IN}$ generated by \MGProx (Algorithm~\ref{aglo:2levelMG_exact}) follows the convergence rate of the proximal gradient method, which is $\mathcal{O}(1/k)$ \cite{beck2017first}.
In Theorem~\ref{thm:rate1k} we show that $\big\{F_0(x^k_0)\big\}_{k \in \IN}$ converges to $F_0^* \coloneqq \inf F_0(x)$ with such a classical rate.

%
%
\begin{theorem}\label{thm:rate1k}
The sequence $\{x_0^k\}_{k \in \IN}$ generated by \MGProx (Algorithm~\ref{aglo:2levelMG_exact}) for solving Problem \eqref{prob:minfg} satisfies
$
F_0(x_0^{k+1}) - F_0^* 
\leq 
\max\big\{8 \delta^2 L_0,\,  F_0(x^1_0) - F_0^* \big\}/k
$,
where $F_0^* = F_0(x^*)$ for $x^* = \argmin F_0$, the point $x^1_0 \in \IRn$ is the initial guess, $L_0$ is the Lipschitz constant of $\nabla f_0$, and $\delta$ is the diameter of the sublevel set $\cL_{\leq F_0(x^1_0)}$ defined in Lemma~\ref{lemma:diamSubLvSet}.
\end{theorem}

Note that we cannot invoke a standard theorem about the convergence of proximal gradient descent such as \cite[Theorem 10.21]{beck2017first}, because we interlace proximal gradient steps with coarse corrections.
Also, we note that all the functions and variables in this subsubsection are at level $\ell=0$ so we omit the subscript.
The constant $L$ should be understood as the Lipschitz constant of $\nabla f(x)$.
The proof is based on standard techniques in first-order methods.
To make the proof more accessible, we divide the proof into four lemmas:
\begin{itemize}
    \item Lemma~\ref{lemma:mgprox_suff_descent}: we derive a sufficient descent inequality for the \MGProx iteration.
    
    \item Lemma~\ref{lemma:quadOverestimator}: we derive a quadratic under-estimator of $F$.
    
    \item Lemma~\ref{lemma:diamSubLvSet}: we give an upper bound for $\| x^k -x^*\|_2$ and $\| y^{k+1} -x^*\|_2$ for all $k$.

    \item Lemma~\ref{lemma:sequence}: we recall a convergence rate for a certain a monotonic sequence. 
\end{itemize}
Using these lemmas, we follow the strategy used in \cite{karimi2017imro} to prove Theorem~\ref{thm:rate1k}.

%
%
\begin{lemma}[Sufficient descent of \MGProx iteration]\label{lemma:mgprox_suff_descent}
For all iterations $k$, we have
\begin{equation}\label{lemma:mgprox_suff_descent_iq}
F(x^{k+1}) - F^*
~\leq~
L\big(\| x^k - x^* \|_2^2 - \| y^{k+1} - x^* \|_2^2\big)/2.
\end{equation}
\end{lemma}
We put the proof in the appendix.
We name the inequality \eqref{lemma:mgprox_suff_descent_iq} sufficient descent because it resembles the sufficient descent property of the proximal gradient iteration \eqref{eqn:suffdescentPGD}. 
Also, by definition, $F(x^{k+1}) \geq F^*$, hence \eqref{lemma:mgprox_suff_descent_iq} implies  $\| x^k - x^* \|_2^2 \geq \| y^{k+1} - x^* \|_2^2$.

%
%
The following lemma is similar to \cite[Lemma 2.3]{beck2009fast} and \cite[Lemma 3, Eq.(5.9)]{karimi2017imro}.
\begin{lemma}[A quadratic function]\label{lemma:quadOverestimator}
For all $x$, we have
\begin{equation}
F(x) - F(x^{k+1})
~\geq~ 
L \langle x^k - y^{k+1}, x - x^k \rangle 
+ L \| y^{k+1} - x^k \|_2^2 / 2.
\label{lemma:QI}
\end{equation}
\end{lemma}
%
%
We put the proof in the appendix.
\begin{lemma}[Diameter of sublevel set]\label{lemma:diamSubLvSet}
At initial guess $x^1 \in \IRn$, define 
\[
\begin{array}{rll}
\cL_{\leq F(x^1)} \hspace{-2mm}
&\coloneqq \big\{ x \in \IRn ~|~ F(x) \leq F(x^1) \big\},
& \textrm{(sublevel set of $x^1$)}
\\
\delta = \text{diam }\cL_{\leq F(x^1)} \hspace{-2mm}
&\coloneqq
\sup \big\{ \| x - y \|_2 ~|~ F(x) \leq F(x^1), F(y) \leq F(y^1) \big\}.
& \textrm{(diameter of $\cL_{\leq F(x^1)}$)}
\end{array}
\]
Then for $x^* \coloneqq \argmin F(x) $, we have
$\| x^k - x^* \|_2 \leq \delta $ and $\| y^k - x^* \|_2 \leq \delta $ for all $k$. 
\begin{proof}
By definition $F(x^*) \leq F(x^1)$.
By the descent property of the coarse correction and proximal gradient updates, we have $F(x^k) \leq F(x^1)$ and $F(y^k) \leq F(x^1)$ for all $k$.
These results mean that $x^k,y^{k+1}$ and $x^*$ are inside $\cL_{\leq F(x^1)}$, therefore both 
$\| x^k - x^* \|_2$ and $\| y^{k+1} - x^* \|_2$ are bounded above by $\delta $.
Lastly, $F$ is strongly convex so $\cL_{\leq F(x^1)}$ is bounded and $\delta < +\infty$.
\end{proof}
\end{lemma}

%
%
\begin{lemma}[Monotone sequence]
\label{lemma:sequence}
For a nonnegative sequence $\{\omega_k\}_{k\in \IN} \rightarrow \omega^*$ that is monotonically decreasing with $\omega_1 - \omega^* \leq 4\mu$ and $\omega_k - \omega_{k+1} \geq (\omega_{k+1} - \omega^*)^2/\mu$, it holds that $\omega_k - \omega^* \leq 4\mu/k$ for all $k$.
%
%
\begin{proof}
By induction.
See proof in \cite[Lemma 4]{karimi2017imro}.
\end{proof}
\end{lemma}

%
%
Now we are ready to prove Theorem~\ref{thm:rate1k}.
\begin{proof}
Rearranging the sufficient descent inequality in  Lemma~\ref{lemma:mgprox_suff_descent} gives
\[
\begin{array}{rcl}
F^* - F(x^{k+1})
&\geq&
L\big( \| y^{k+1} - x^* \|_2^2  -    \| x^k - x^* \|_2^2\big)/2
\\
&=&
L \big(\| y^{k+1} - x^* \|_2  - \| x^k - x^* \|_2 \big)
\big( \| y^{k+1} - x^* \|_2  + \| x^k - x^* \|_2 \big)/2
\\
&\geq&
-L \| x^k - y^{k+1} \|_2 \big(\| x^k - x^* \|_2 + \| y^{k+1} - x^* \|_2\big)/2,
\end{array}
\]
where the last inequality is by the triangle inequality
$\| y^{k+1} - x^* \|_2  + \| x^k -  y^{k+1} \|_2 \geq   \|  x^* - x^k \|_2$.
Rearranging the inequality gives
\begin{equation}
\| x^k - y^{k+1} \|_2 
~\geq~
\dfrac{-2}{L}
\frac{F^* - F(x^{k+1}) }{\| x^k - x^* \|_2 + \| y^{k+1} - x^* \|_2}
=
\dfrac{2}{L}
\frac{F(x^{k+1}) - F^* }{\| x^k - x^* \|_2 + \| y^{k+1} - x^* \|_2}.
\label{thm:rate_4}
\end{equation}
Applying Lemma~\ref{lemma:diamSubLvSet} to 
\eqref{thm:rate_4} gives
\begin{equation}
\| x^k - y^{k+1} \|_2 
\geq
\big(F(x^{k+1}) - F^*\big)/(\delta L).
\label{thm:rate_5}
\end{equation}
Note that \eqref{thm:rate_5} implies that if the fine sequence converges ($x^k = y^{k+1}$), then $F(x^{k+1}) = F^*$.

Applying Lemma~\ref{lemma:quadOverestimator} with $x = x^k$ gives 
$
F(x^k) - F(x^{k+1}) 
\geq
L \| y^{k+1} - x^k \|_2^2/2
\overset{\eqref{thm:rate_5}}{\geq}
(F(x^{k+1}) - F^*)^2 / (2\delta^2 L)$.
This inequality shows that the sequence $\{ \omega_k \}_{k\in \IN}$ with $\omega_k \coloneqq F(x^k)$ satisfies the condition $\omega_k - \omega_{k+1} \geq (\omega_k - \omega^* )^2/\mu$ in Lemma~\ref{lemma:sequence}.
To complete the proof, applying Lemma~\ref{lemma:sequence} to the monotonically decreasing sequence $\{F(x^k)\}_{k \in \IN}$ with $\mu = 2\delta^2L$ gives
$F_0(x^{k+1}_0) - F^*_0 
\leq
\max\big\{\, 8 \delta^2 L_0,  F_0(x^1_0) - F_0^* \,\big\} /k
$, where we put back the subscript 0.
\end{proof}

Theorem~\ref{thm:rate1k} shows that $\{F_0(x^k_0)\}_{k \in \IN}$ for solving Problem \eqref{prob:minfg} satisfies a sublinear convergence bound of $\cO(1/k)$.
Below we show that $\{F_0(x^k_0)\}_{k \in \IN}$ satisfies a linear convergence bound.

%
%
\subsubsection{Linear convergence rate}
All the functions and variables here are at level 0 so we omit the subscripts.
Now we show that $\big\{F( x^k) \big\}_{k \in \IN}$ converges to $F^*$ with a linear rate using the \emph{Proximal Polyak-{\L}ojasiewicz inequality} \cite[Section 4]{karimi2016linear}. 
The function $F$ in Problem \eqref{prob:minfg} is called ProxPL, if there exists $\mu > 0$, called the ProxPL constant, such that 
\begin{equation}
\cD_g(x,L)
~\geq~
2 \mu \big( F(x) - F^* \big)  \qquad  \forall x,
\tag{ProxPL}
\label{def:ProxPL}
\end{equation}
\begin{equation}
\cD_g(x,\alpha) 
\coloneqq
-2\alpha \min_z \,
\Big\{
{\textstyle\frac{\alpha}{2}} \| z - x \|_2^2 
+ \big\langle z - x, \nabla f(x) \big\rangle
+ g(z) - g(x)
\Big\}
.
\label{def:dg}
\end{equation}
Intuitively, $\cD_g$ is defined based on the proximal gradient operator:
\[
\prox_{\frac{1}{L}g} \Big( x - \frac{\nabla f(x)}{L} \Big)
=
\underset{z}{\argmin} \, \frac{L}{2} \| z -  x \|_2^2
+ \big\langle z - x, \nabla f(x) \big\rangle
+ g(z) - g(x).
\]
It has been shown in \cite{karimi2016linear} that if $f$ in \eqref{prob:minfg} is $\mu$-strongly convex, then $F$ is $\mu$-ProxPL.
Now we prove the linear convergence rate of Algorithm~\ref{aglo:2levelMG_exact}.
Note that a standard result such as \cite[Theorem 10.29]{beck2017first} on convergence of proximal gradient for strongly convex functions is not directly applicable because, as mentioned above, we interleave proximal gradient steps with coarse correction steps.
%
%
\begin{theorem}\label{thm:mgprox_converge_proxPL}
Let $x^1_0$ be the initial guess of the algorithm,
$x_0^* = \argmin \, F_0(x)$ and $F^*_0 \coloneqq F_0(x_0^*)$.
The sequence $\{x_0^k\}_{k \in \IN}$ generated by \MGProx (Algorithm~\ref{aglo:2levelMG_exact}) for solving Problem \eqref{prob:minfg} satisfies
$
F_0(x_0^{k+1}) - F_0^* 
\leq 
\Big( 1-\frac{\mu_0}{L_0} \Big)^k \big( F_0(x_0^1) - F_0^* \big).
$

\begin{proof}
First, by assumption $F$ is strongly convex, so $F$ is $\mu$-ProxP\L\, with $\mu > 0$ and
\[
\renewcommand*{\arraystretch}{1.33}
\begin{array}{rcl}
F(x^{k+1})
&\overset{\eqref{eqn:suffdescentPGD}}{\leq}& \hspace{-2mm}
F(z^{k+1})
~\overset{\eqref{iq:coarse_descent_on_F}}{\leq}~ 
F(y^{k+1})  
~=~
f(y^{k+1}) + g(y^{k+1}) + g(x^k) - g(x^k)
\\
&\overset{\eqref{def:Lsmooth}}{\leq}&  \hspace{-2mm}
f(x^k) +  \langle  \nabla f(x^k), y^{k+1} - x^k \rangle + \frac{L}{2} \|  y^{k+1} - x^k \|_2^2 + g(x^k) +  g(y^{k+1})  -  g(x^k) 	
\\
&=& \hspace{-2mm}
F(x^k) +  \langle  \nabla f(x^k), y^{k+1} - x^k \rangle + \frac{L}{2} \|  y^{k+1} - x^k \|_2^2 +  g(y^{k+1})  -  g(x^k) 
\\
&=& \hspace{-2mm}
F(x^k) - \frac{1}{2L} \underbrace{(-2L)\Big( 
\langle  \nabla f(x^k), y^{k+1} - x^k \rangle 
+ {\textstyle \frac{L}{2}} \|  y^{k+1} - x^k \|_2^2 +  g(y^{k+1})  -  g(x^k) 
\Big)}_{\overset{\eqref{def:dg}}{ = } \, \cD_g(x^k, L), \, \text{ since }y^{k+1} \,\coloneqq\, \underset{z}{\argmin} \, \frac{L}{2} \| z -  x^k \|_2^2
+ \langle z - x^k, \nabla f(x^k)\rangle
+ g(z) - g(x^k).}
\\
&\overset{\eqref{def:ProxPL}}{\leq}& \hspace{-2mm}
F(x^k) - \frac{\mu}{L} \big(F(x^k) - F^*\big).
\end{array}
\]
Adding $-F^*$ on both sides of the inequality gives
$ 
F(x^{k+1}) - F^* 
\leq \big( 1 - \mu/L\big) \big(F(x^k) - F^*\big)
$.
Applying this inequality recursively completes the proof.
\end{proof}
\end{theorem}

\begin{remark}
We now give several remarks about the result.
\begin{itemize}
\item Convergence rate: for a $\mu_0$-strongly convex and $L_0$-smooth $f_0$, we have $0 < \mu_0 \leq L_0$ and $0 \leq 1 - \mu_0/L_0 < 1$.
Depending on the value of $\mu_0$, for $k$ not too large, the sublinear convergence rate $1/k$ from Theorem~\ref{thm:rate1k} gives a better bound than the linear rate $(1-\mu_0/L_0)^k$ from Theorem~\ref{thm:mgprox_converge_proxPL}.
This is the case when $\mu_0 \ll L_0$.

\item Since $x^*$ is unique, we also conclude that the sequence $\{x^k_0\}_{k \in \IN}$ converges to $x^*$. 
\end{itemize}
\end{remark}
Lastly, we emphasize that the bounds in Theorem~\ref{thm:rate1k} and Theorem~\ref{thm:mgprox_converge_proxPL} are loose bounds for Algorithm~\ref{aglo:2levelMG_exact} since we only show that the coarse correction can guarantee the descent condition $F_0(z_0^{k+1}) \leq F_0(y_0^{k+1})$ but not a stronger sufficient descent condition.

%
%
\subsubsection{On the selection of \texorpdfstring{$\tau$}{}}\label{sec:mgprox:subsec:tau}
Recall that the tau vector comes from a set:
\[
\tau_{0 \rightarrow 1}^{k+1} 
~\overset{\eqref{def:tau_mgprox2}}{\in}~
\underline{\tau_{0 \rightarrow 1}^{k+1}}
\coloneqq 
\underbrace{ \nabla f_1(y_1^{k+1}) - R \partial F_0(y_0^{k+1}) }_{\textrm{singleton}} 
+
\underbrace{ \partial g_1(y_1^{k+1} ) }_{\textrm{set-valued}}.
\]
We emphasize that for our theoretical results to hold we can choose any value in the set $\underline{\tau_{0 \rightarrow 1}^{k+1}}$ to define $\tau_{0 \rightarrow 1}^{k+1}$ in Algorithm~\ref{aglo:2levelMG_exact}.
First, the results of in Theorems~\ref{thm:fpp} and~\ref{thm:coarse_dir_descent} hold for any $\tau$ in the set.
Second, all the convergence bounds (Theorems~\ref{thm:rate1k} and \ref{thm:mgprox_converge_proxPL}) only contain constants at level $\ell = 0$ and are independent of the choice of tau vector.

\paragraph{Optimal tau selection seems difficult}
Recall the two steps in the algorithm related to the coarse correction,
\[
\begin{array}{rcl}
x_1^{k+1}(\tau)
&=& \underset{\xi}{\argmin}\, F_1(\xi) - 
\Big\langle
\underbrace{\nabla f_1(y_1^{k+1}) - R \partial F_0(y_0^{k+1})  + \partial g_1(y_1^{k+1}) }_{\ni \tau}
,\,\xi 
\Big\rangle,
\\
x_0^{k+1}(\tau) 
&=& 
y_0^{k+1} + \alpha(\tau) P \big(      x_1^{k+1}(\tau) - y_1^{k+1}       \big),
\end{array}
\]
where $x_0^{k+1}, x_1^{k+1}$ and $\alpha$ are all a function of $\tau$.
Now it seems tempting to ``optimally tune'' $\tau_{0 \rightarrow 1}^{k+1}$ so that it maximizes the gap $F_0(y_0^{k+1}) - F_0(x_0^{k+1}) $:
\[
\tau_{0\rightarrow 1}^{k+1}
\in \underset{\tau \,\in\, \underline{\tau_{0 \rightarrow 1}^{k+1}}}{\argmax} \, F(y_0^{k+1}) - F_0\big(x_0^{k+1}(\tau) \big)
= \underset{\tau}{\argmin} \,  F_0\big(x_0^{k+1}(\tau) \big).
\]
However, this problem generally has no closed-form solution and it is intractable to solve numerically.
In the experiments we verified that the sequence produced by \MGProx converges for different values of $\tau$ confirming the theory. 

\section{FastMGProx: MGProx with Nesterov's acceleration}\label{sec:NestMGProx}
In this section we show that, by treating the \MGProx update as a single iteration, we can embed this update within Nesterov's estimating sequence framework to derive an accelerated method with the optimal $\cO(1/k^2)$ convergence.
In this section we are focusing on level $\ell=0$ and we hide some of the subscripts $0$ for clarity.

Algorithm~\ref{aglo:fast_2levelMG_exact} shows the Nesterov's accelerated \MGProx that we call \texttt{FastMGProx}.
First we introduce some compact notation.
We denote $\Big(\text{MGProx}\circ\prox\Big)(w)$ as applying a proximal gradient update on $F_0$ at the point $w$ and then followed by a MGProx update process (i.e., steps (i)-(vi) in Algorithm~\ref{aglo:2levelMG_exact}).
The algorithm introduces two scalar sequences $\alpha^k, \gamma^k$ and two auxiliary vector sequences $y^k, z^k$.
It is important not to confuse these symbols with those presented in Section~\ref{sec:mgprox}.

\begin{algorithm}
\caption{Fast\MGProx for an approximate solution of \eqref{prob:minfg}}
\label{aglo:fast_2levelMG_exact}
\begin{algorithmic}
\STATE{Initialize $R$ and $P$, $z^0 = x^0$ (auxiliary sequence)$, \gamma^0 > 0$ (extrapolation parameter)}
\FOR{$k = 1,2,\dots$}
\STATE{(i)~~~~
Solve $\alpha^k \in ~]0,1[$ from $L_0 (\alpha^k)^2 = (1-\alpha^k) \gamma^k $ 
,
$\gamma^{k+1} = (1-\alpha^k)\gamma^k$
\hfill parameter
}
\STATE{(ii)~~~~
$y^k ~~~= \alpha^k z^k +  (1-\alpha^k)x^k$
\hfill Nesterov's extrapolation
}
\STATE{(iii)~~~~ 
$x^{k+1} = \big(\text{MGProx} \circ \prox\big)\big( y^k - \frac{1}{L_0}\nabla f_0(y^k) \big)$
\hfill prox-grad step with MGProx
}
\STATE{(iv)~~~~
$z^{k+1} = z^k - \frac{\alpha^k}{\gamma^{k+1}} \frac{y^k - x^{k+1}}{L_0}$
\hfill updating the auxiliary sequence
}
\ENDFOR
\end{algorithmic}
\end{algorithm}

By the sufficient descent lemma of proximal gradient update 
(Lemma~\ref{lemma:suff_descent_proxgrad})
and the descent lemma of MGProx update 
(inequality \eqref{iq:coarse_descent_on_F}), we can guarantee that for $x^k, y^k$ in Algorithm~\ref{aglo:fast_2levelMG_exact} we have $ F(x^{k+1}) \leq F(y^k) - \| G(y^k) \|_2^2/(2L)$,
where $G(y^k)$ is the proximal gradient map of $F$ at $y^k$, see  \eqref{eqn:suffdescentPGD}.
This inequality forms the basis of the \texttt{FastMGProx}: using Nesterov's framework of estimating sequence \cite{nesterov2003introductory}, we have the following theorem

\begin{theorem}\label{thm:FMGProx_fastr_ate}
For the sequence $\{x^k\}$ produced by Algorithm~\ref{aglo:fast_2levelMG_exact}, we have the convergence rate $F(x^k) - F^* \leq \cO(1/k^2)$.
To be exact,
\[
F(x^k) - F^*   ~\leq~ 
\dfrac{4L_0
\Big(  F(x^0) - F^* + \dfrac{\gamma^0}{2}\| x^0 - x^* \|_2^2 \Big)
}{
\Big(2\sqrt{L_0} - \sqrt{\gamma^0}\Big)^2
+ 2\gamma^0\Big( 2\sqrt{L_0}-1 \Big)k
+ \gamma^0 k^2
}.
\]
\end{theorem}

We present the proof of Theorem~\ref{thm:FMGProx_fastr_ate} in Appendix \ref{sec:appendix}.
Below we give some key points about Algorithm~\ref{aglo:fast_2levelMG_exact}.
First, the convergence rate $F(x^k) - F^* \leq \cO(1/k^2)$ is optimal in the function-gradient model of \cite{nesterov2003introductory}.
Next, the work \cite{parpas2017multilevel} also achieved the same optimal rate.
However, inspecting the algorithm of \cite{parpas2017multilevel}, it relies on two things:
\begin{itemize}
    \item It bypassed all the technical challenges caused by the non-smoothness of the coarse problem by using a smooth approximation.
    In contrast, \MGProx and \texttt{FastMGProx} do not use smoothing.


    \item There is a safe-guarding IF statement in the algorithm.
    It means that when the combined effect of Nesterov's extrapolation and multigrid process produces a bad iterate, that whole iteration will be discarded and replaced by a simple proximal gradient iteration.
    This gives a sign of possible bad interaction between  Nesterov's acceleration and multigrid.
    In contrast, \MGProx and \texttt{FastMGProx} have no such IF statement.
\end{itemize}

\section{A multi-level MGProx}\label{sec:multilevel}
Now we generalize the 2-level \MGProx to multiple levels.
The 2-level \MGProx method constructs a coarse problem at level ($\ell = 1$), and uses the solution of such problem to help solve the original fine-level problem ($\ell = 0$).
If the fine problem has a large problem size, solving the coarse problem exactly is generally expensive.
Hence it is natural to consider applying multigrid recursively until the coarse problem on the coarsest level is no longer expensive to solve. 
An $L$-level \MGProx cycle with a V-cycle structure is shown in Algorithm~\ref{algo:LlevelMG_exact}.
We clarify the naming of the variables in the algorithm as follows: at each iteration $k$, we have
$x_\ell^k$: variable before pre-smoothing on level $\ell$;
$y_\ell^{k+1}$: variable after pre-smoothing on level $\ell$;
$z_\ell^{k+1}$: variable after coarse-grid correction on level $\ell$; and
$w_\ell^{k+1}$: variable after post-smoothing on level $\ell$.
Note that, to obtain a well-defined recursion in Algorithm~\ref{algo:LlevelMG_exact}, we choose the superscript for the $x$ variables equal to $k$ on all levels. In the 2-level algorithm we chose a different convention, writing the $x$ variable on level 1 as $x_{1}^{k+1}$.

\begin{algorithm}
\caption{$L$-level \MGProx with V-cycle structure for an approximate solution of \eqref{prob:minfg}}
\label{algo:LlevelMG_exact}
\begin{algorithmic}
\STATE{Initialize $x^1_0$ and the full version of $R_{\ell \rightarrow \ell+1}, P_{\ell+1 \rightarrow \ell }$ for 
$\ell \in \{0,1,\dots,L-1\}$
}

\FOR{$k = 1,2,\dots$}
    \STATE{Set $\tau_{-1 \rightarrow 0 }^{k+1}=0$}
    \FOR{$\ell = 0,1,\dots, L-1$}
            \STATE{$y^{k+1}_\ell \,~~~= \prox_{\frac{1}{L_{\ell}}g_{\ell}} 
             \bigg (
             x^{k}_\ell - \dfrac{\nabla f_\ell(x^k_\ell)  - \tau_{\ell-1 \rightarrow \ell }^{k+1}}{L_\ell}
             \bigg)$ \hfill pre-smoothing
             }         
             \STATE{$x^{k}_{\ell+1} \,~~~= R_{\ell \rightarrow \ell+1}(y^{k+1}_{\ell}) \, y^{k+1}_{\ell} $
             \hfill restriction to next level \vspace{0.5mm}
             }
             \STATE{
            $\tau_{\ell \rightarrow \ell+1}^{k+1} \in \underline{\partial F_{\ell+1}(x^{k}_{\ell+1} )}  - R_{\ell \rightarrow \ell+1}(y^{k+1}_{\ell}) \, \underline{\partial  F_\ell(y^{k+1}_{\ell})}$
            \hfill create tau vector ~
            }
    \ENDFOR
    
             \STATE{$w_L^{k+1} = \underset{\xi}{\argmin} \, \Big\{\, F_L^{\tau}(\xi) \coloneqq F_L(\xi) - \langle \tau_{L-1 \rightarrow L}^{k+1}, \xi \rangle \,\Big\}$
             \hfill solve the level-$L$ coarse problem}

    \FOR{$\ell = L-1, L - 2, \dots, 0$}
        \STATE{
        $z_\ell^{k+1} = y^{k+1}_\ell  + \alpha  P_{\ell+1 \rightarrow \ell}\big(  w^{k+1}_{\ell+1}  - x^{k}_{\ell+1} \big)$
        \hfill coarse correction
        }
        \STATE{
        $w^{k+1}_\ell = \prox_{\frac{1}{L_\ell}g_\ell} 
        \bigg(
        z^{k+1}_\ell - \dfrac{\nabla f_\ell(z^{k+1}_\ell) -\tau_{\ell-1 \rightarrow \ell }^{k+1}}{L_\ell} 
        \bigg)$
        \hfill post-smoothing
        }
    \ENDFOR
    \STATE{$x^{k+1}_{0} \,~~~= w^{k+1}_0$ \hfill update the fine variable}
\ENDFOR
\end{algorithmic}
\end{algorithm}

Here are remarks about Algorithm~\ref{algo:LlevelMG_exact}.
First, $L_\ell$ is the Lipschitz constant of $\nabla f_\ell$.
Then,
\begin{itemize}
\item At level $\ell \neq L $, we are essentially performing two proximal gradient iterations (pre-smoothing + post-smoothing) and a coarse correction.
    At the coarsest level $\ell = L$, we perform an exact update by solving the coarse problem exactly.
    
\item From the traditional wisdom of classical multigrid, more than one pre-smoothing and post-smoothing steps can be beneficial to accelerate the overall convergence.
We implemented such multiple smoothing steps in the numerical tests.
\end{itemize}

\begin{remark}[Convergence of Algorithm~\ref{algo:LlevelMG_exact}]
Regarding the finest level function value $\{F_0(x_0^k)\}_{k \in \IN}$, Theorem~\ref{thm:rate1k} and Theorems~\ref{thm:mgprox_converge_proxPL} and \ref{thm:mgprox_converge_proxPL} all hold for the multilevel Algorithm~\ref{algo:LlevelMG_exact}, since the angle condition of the coarse correction (Theorem~\ref{thm:coarse_dir_descent}) also holds for multilevel \MGProx when the coarse problem is solved inexactly.
To be specific, the last inequality $F_1^{\tau}(x^{k+1}_1) - F_1^{\tau}(y^{k+1}_1) < 0$ in the proof of Theorem~\ref{thm:coarse_dir_descent} holds when the coarse problem on $x^{k+1}_1$ is solved inexactly by a combination of proximal gradient iterations and a coarse-grid correction with line search.
\end{remark}

%
%
\paragraph{Approximate per-iteration computational complexity}\label{sec:multilevel:cost}
The cost of one prox-grad update on the finest level scales with $n_0$. 
The cost of all prox-grad operations in one iteration of a $L$-level \MGProx is 
$
2 ( 1 + r + \dots + r^{L-1} ) n_0 = 2\frac{1-r^L}{1-r} n_0
$,
where 
$r$ is a reduction factor.
In the experiment $r=\frac{1}{4}$, and all the prox-grad steps per V-cycle iteration amount to at most 
$
 \frac{8}{3}(1-\frac{1}{4^L})
\leq 
2.67
$
 times the cost of performing one fine prox-grad update.

\paragraph{Multi-level FastMGProx}
Similar to the $L$-level \texttt{MGProx}, we can also propose a $L$-level \texttt{FastMGProx}, which is not shown here.

\section{Numerical results}\label{sec:exp}
We test \MGProx on an Elastic Obstacle Problem (EOP), to be reviewed below, then we present the test results in \cref{sec:exp:subsec:result}.
%
%
\subsection{Elastic Obstacle Problem (EOP)}
The EOP \cite{brandt1983multigrid,mandel1984multilevel} was motivated by physics \cite{rodrigues1987obstacle}, see \cite[Section 4]{tran20151} for related problems.
Here we solely focus on solving EOP as a nonsmooth optimization problem in the form of \eqref{prob:minfg} and use it to demonstrate the capability of \texttt{MGProx}.

As elastic potential energy is proportional to the surface area, we have the problem
\begin{equation}
 \min_{ u } \iint_\Omega \sqrt{1 +\| \,\nabla u \, \|^2_{L^2}} dxdy
+ \lambda \iint_{\Omega} \| \,(\phi - u)_+ \, \|_{L^1} dxdy
~~~\text{s.t.}~~~ u = 0 \text{ on } \partial \Omega,
\tag{continuum EOP}
\label{prob:EOP}
\end{equation}
where $\nabla u : \Omega \rightarrow \IR^2$ is the gradient field of $u$,
the norms $\| \cdot \|_{L^2}$ and $\| \cdot \|_{L^1}$ are the $L^2$ and $L^1$ norm for functions, resp..
We take $\phi \leq 0$ on $\partial \Omega$, so the boundary condition $u = 0$ on $\partial \Omega$ makes sense.
Here the non-penetration constraint is represented by the penalty term, where 
$(\cdot)_+ \coloneqq \max\{\cdot, 0\}$ is taken element-wise,  and $\lambda >0$ is a penalty parameter.

\paragraph{Most multigrid algorithms solve an approximated version of EOP}
Many methods have been proposed to solve EOP: adaptive finite element methods \cite{hoppe1994adaptive},
penalty methods \cite{scholz1984numerical, tran20151} and level set methods \cite{majava2004level}.
For the multigrid methods mentioned in the introduction,
 many of them do not solve \eqref{prob:EOP} efficiently. 
\begin{itemize}
\item They only solve an approximate version of EOP: 
near $x=0$ the function $\sqrt{1+x^2}$ has a Taylor series $1 + x^2/2 + o(x^3)$, ignoring the 1 and higher order terms, such linearization replace the integral of \eqref{prob:EOP} by $\frac{1}{2}  \int_\Omega \| \nabla u  \|^2_{L^2} \textrm{dx}$.

\item They only solve the box-constrained version of EOP, i.e., they are not designed to handle the nonsmooth $L^1$ penalty.
In fact \cite{kocvara2016first} pointed out it is general hard to extend multigrid from dealing with box constraint to general nonsmooth functions.

\item They are not FAS but Newton-MG (see \cref{sec:intro:subsec:review}).
\end{itemize}

%
%
\subsection{The optimization problem}
We now discuss how to derive the optimization problem by discretizing the integral in \eqref{prob:EOP} to obtain a problem in the form of \eqref{prob:minfg}.
On a grid of $N\times N$ internal points with $\Delta x=\Delta y=h= \frac{1}{N+1}$ on $\Omega$, we let $\U \in \IR^{N \times N}$ with $\U(i,j) = u( ih, jh)$ with $i,j$ ranging from $1$ to $N$.
Here small italic symbol $u(x,y)$ denotes the continuum variable in the infinite dimensional space, and the capital bold symbol $\U$ denotes a $N$-by-$N$ matrix obtained by finite discretization of $u$.
Let $\vec$ denotes vectorization, 
then $\u \coloneqq \vec(\U) \in \IR^{N^2}$ is the optimization variable.
Discretizing the integral in \eqref{prob:EOP} gives
\begin{equation}\label{prob:EOP-0}
\argmin_{\u \in \IR^{N^2}} 
h^2 \sum_{i=1}^N \sum_{j=1}^N \sqrt{
1
+\big( \D_{(i,j),: } \u \big)^2
+\big( \E_{(i,j),: } \u \big)^2
}
+ 
h^2 \lambda \| (\bphi - \u)_+ \|_1
\tag{EOP-0}
\end{equation}
where the two matrices 
$\D \in \IR^{N^2 \times N^2}, \E \in \IR^{N^2 \times N^2}$ are the first-order forward difference operators that approximate $\partial/\partial x$ and $\partial/\partial y$ respectively, defined using the coordinate index $(i,j)$ of $\D$ and the coordinate index $(k,l)$ of $\U$ as
\[
\D_{(i,j),(k,l)} 
=
\begin{cases}
1/h & i=k\mbox{ and } j =l+1, 
\\
-1/h &  i=k\mbox{ and } j =l, 
\\
0            & \mbox{else},
\end{cases}
~~
\E_{(i,j),(k,l)}
=
\begin{cases}
1/h  & i=k+1\mbox{ and } j =l, 
\\
-1/h & i=k ~~~~~~\mbox{ and } j =l,
\\
0            & \mbox{else}.
\end{cases}
\]
The notation $\D_{(i,j),: }$ refers to the $(i,j)$-th row of $\D$.
The vector $\bphi \in \IR^{N^2}$ is the discretization of $\phi$ in \eqref{prob:EOP}, and $\| \cdot \|_1$ is the $\ell_1$-norm.

Define $\F \in \IR^{2N^2 \times N^2}$ by stacking rows of $\D$ and $\E$ together as in \eqref{def:psi_F}.
Ignoring the constant factor $h^2$ in \eqref{prob:EOP-1} 
yields the problem
\begin{equation}\label{prob:EOP-1}
\argmin_{\u \in \IR^{N^2}} 
f_0(\u) + g_0(\u)
\coloneqq
\sum_{i=1}^N \sum_{j=1}^N \psi\Big(
\F_{(i,j),:} \u 
\Big)
+ 
\lambda \| (\bphi - \u)_+ \|_1,
\tag{EOP-1}
\end{equation}
where
\begin{equation}\label{def:psi_F}
\psi : \IR^2 \rightarrow \IR : (s,t) \mapsto \sqrt{1+s^2+t^2}
,
\qquad
\F_{(i,j),:} \coloneqq 
\begin{bmatrix}
\D_{(i,j),:}
\\
\E_{(i,j),:}
\end{bmatrix}.
\end{equation}

To avoid confusing the notion of gradient in the finite dimensional Euclidean space with the gradient operator $\nabla $ in \eqref{prob:EOP}, from now on we denote $\grad f_0$ the gradient of $f_0$ in \eqref{prob:EOP-1}.
We denote $\Hess f_0$ the Hessian of $f_0$.
By chain rule and the fact that
\begin{equation}\label{grad_psi}
\grad \, \psi (s,t) = 
\left.
\begin{bmatrix}
\dfrac{\partial}{\partial x} \psi (x,y)\vspace{1mm}
\\
\dfrac{\partial}{\partial y} \psi (x,y)
\end{bmatrix} \, \right|_{x=s, y =t}
=
\begin{bmatrix}
\dfrac{s}{\sqrt{1+s^2+t^2}}\vspace{1mm}
\\
\dfrac{t}{\sqrt{1+s^2+t^2}} 
\end{bmatrix}
\in \IR^2,
\end{equation}
so 
$
\grad f_0(\u) 
= \displaystyle 
\sum_{i,j} 
\underbrace{
\F_{(i,j),: } ^\top
}_{N^2 \text{-by-} 2}
\underbrace{
\grad\,\psi \Big( \F_{(i,j),: } \u \Big)
}_{2 \text{-by-} 1}
=
\sum_{i,j} 
\begin{bmatrix}
\D_{(i,j),: } \vspace{0.3mm}
\\
\E_{(i,j),: }
\end{bmatrix}^\top
\grad\, \psi \Big( \D_{(i,j),: } \u,   ~ \E_{(i,j),: }\u \Big)
\in \IR^{N^2}
$ is
\begin{equation}\label{grad_f}
\grad f_0(\u) 
~\overset{\eqref{def:psi_F}}{=} ~
\sum_{i=1}^N 
\sum_{j=1}^N 
\tilde{\psi}_{ij}^D
\D_{(i,j),: } ^\top
+
\tilde{\psi}_{ij}^E
\E_{(i,j),: }^\top
~=~ 
\D^\top \bar{\bpsi}^D
+
\E^\top \bar{\bpsi}^E
\end{equation}
where 
$ \displaystyle
\tilde{\psi}_{ij}^D 
\coloneqq 
\dfrac{  [\D\u]_{ij} }{ 
\sqrt{1+ \big( [\D\u]_{ij} \big)^2 + \big( [\E\u]_{ij} \big)^2}
}
$ 
and 
$
\tilde{\psi}_{ij}^E 
\coloneqq 
\dfrac{ [\E\u]_{ij} }{ 
\sqrt{1+ \big( [\D\u]_{ij} \big)^2 + \big( [\E\u]_{ij} \big)^2}
}
$.

By proposition~\ref{prop:strcvx},  $\grad f_0$ is Lipschitz and problem \eqref{prob:EOP-1} is strongly convex and thus: 
1. it has a unique minimizer, and 
2. it is within the framework of Problem \eqref{prob:minfg} so \MGProx can be used.

\begin{proposition}\label{prop:strcvx}
For $\u$ on a bounded domain, the function $f_0$ in \eqref{prob:EOP-1} is $L_0$-smooth (i.e., the gradient $\grad f_0$ is $L_0$-Lipschitz) and strongly convex, with $L_0$ provided below.
\end{proposition}
\begin{proof}
Consider $\psi(s,t)$ defined by  \eqref{def:psi_F}, $\grad\, \psi(t)$ is given by \eqref{grad_psi} and the Hessian is
\[
\Hess\, \psi(s,t)= \dfrac{1}{(1+s^2+t^2)^{3/2}}
\left(
\begin{array}{cc}
  1+t^2 & -st \\
  -st & 1+s^2
\end{array}
\right).
\]
This shows that $\psi$ is 1-smooth and strongly convex on any bounded domain.

Now recall the function $f_0(\u)$ can be written as $f_0(\u) = \sum_{i,j} \psi\Big(\F_{(i,j),:} \u \Big)$
as in \eqref{prob:EOP-1} with $\F_{(i,j),:}$ defined by \eqref{def:psi_F}.
Since a composition of a convex function and linear function is convex, and the sum of convex functions is convex, this shows $f_0$ is convex.
Furthermore,
\[
\Hess\, f_0(\u)
~=~
\sum_{i,j=1}^n \Hess\,\psi\Big( \F_{(i,j),:} \u\Big)
\cdot \F_{(i,j),:}^\top \F_{(i,j),:}.
\]
Then we compute that the Lipschitz constant $L_0$ of $\grad f_0$ is at most $\displaystyle\sum_{i,j} \|  \F_{(i,j),:}^\top \F_{(i,j),:} \|_2$.
Since $\|  \F_{(i,j),:} \|_2 = \sqrt{3}/h$, we can take $L_0=\sqrt{3}N^2/h$.

For strong convexity, we see that the Hessian is a weighted sum of rank-2 matrices with positive weights.
Suppose $\z$ lies in the null space of all the rank-2 matrices.  
In this case, $\z$ must be identically $\mathbf{0}$ as the following argument shows.  
If for some $(i,j)$ it holds that $\D_{(i,j),:}\z=\mathbf{0}$, then $z_{i+1,j}=z_{i,j}$.
So, if $\z$ is in the null space of all the $\D_{(i,j),:}$'s, then $\z$ is constant in the $x$-direction.  
Similarly, if $\z$ is in the null space of all the $\E_{(i,j),:}$'s, then $\z$ is constant in the $y$-direction.
The operators $\D,\E$ applied at the boundaries yield that $\z$ is identically zero.

Thus, $\Hess f_0(\u)$ is positive definite for all $\u$.
Since the Hessian
is continuous, this means that its smallest eigenvalue has a positive lower
bound over any bounded domain, i.e., $f_0$ is strongly convex on any bounded
domain.
\end{proof}

\begin{remark}\label{remark:Kmu_unknown}
It is expensive to evaluate $f_0$ and $\grad f_0$ for \eqref{prob:EOP-1}, since they are sum of $N^2$ nonlinear terms involving vectors in $\IR^{N^2}$.
Furthermore, based on Proposition~\ref{prop:strcvx},
\begin{itemize}
    \item a tight closed-form expression of the global Lipschitz constant of $\grad f_0$ is unknown to us,
    \item $f_0$ is $\mu_0$-strongly convex but the strong convexity parameter $\mu_0$ is unknown to us.
\end{itemize}

\end{remark}

\subsubsection{MGProx on EOP}\label{sec:exp:subsec:mgprox_eop}
The subdifferential and proximal operator of $\lambda \|(\c-\v)_+ \|_1$ are
\begin{equation}
\Big[\partial \| (\c-\v)_+ \|_1 \Big]_i
= \begin{cases}
-1 & v_i < c_i 
\\
[-1,0]  & v_i = c_i
\\
0 & v_i > c_i
\end{cases}
,
~~
\Big[\prox_{\lambda \| (\c-\,\cdot\,)_+ \|_1}(\v)\Big]_i
= 
\begin{cases}
v_i + \lambda  &  v_i + \lambda < c_i
\\
c_i  &  v_i \leq c_i \leq v_i + \lambda 
\\
v_i &  v_i > c_i 
\end{cases}.
\label{eqn:subdiff_prox_maxpen}
\end{equation}
For a matrix $\U$ at a resolution level $\ell$, the (full) restriction of $\U$, denoted as $\hat{\U} \coloneqq \overline{\cR}(\U)$, can be defined by a full weighting operator with the kernel defined by 
$
\frac{1}{8} [1\, 2\, 1]^\top \otimes [1\, 2\, 1]
$
where $\otimes$ is tensor product.
Using a block-tridiagonal matrix $\overline{\R}$ and vectorization $\vec$, the expression $\hat{\U} \coloneqq \overline{\cR}(\U)$ can be written as $\hat{\u} \coloneqq \vec(\hat{\U}) = \overline{\R}\u \coloneqq \overline{\R}\vec(\U)$.
For the adaptive version of $\R$, we follow Definition~\ref{def:adaptiveR}.
For the prolongation matrix $\P \coloneqq c \R^\top$, we take $c=2$.


%
%
\subsection{Test results}\label{sec:exp:subsec:result}
We now present the test results.

\paragraph{Experimental setup}
We take
$\phi(x,y) = \max\{0,\sin(x)\} \max\{0,\sin(y)\}$ where $x,y \in [0,3\pi]$.
We initialize $x_0^1$ as a random nonnegative vector and compute the initial function value $F_0(x_0^1)$ and the initial norm of the proximal gradient map $\|G_0(x^ \textrm{ini}_0)\|_2$.
We stop the algorithm using the proximal first-order optimality condition: we stop the algorithm if 
$\|G_0(x^k_0)\|_2 / \|G_0(x^\textrm{ini}_0)\|_2 \leq 10^{-16}$.
The experiments are conducted in MATLAB R2023a on a MacBook Pro (M2 2022) running macOS 14.0 with 16GB memory\footnote{The code is available at \url{https://angms.science/research.html}}.
We report the value $ ( F_0(x^k_0)  - F_0^\textrm{min} )/  F_0(x^\textrm{ini}_0) $, where $F_0^\textrm{min}$ is the lowest objective function value achieved among all the tested methods.

\paragraph{MGProx setup}
We run a multilevel \MGProx with $N_s$ number of pre-smoothing and post-smoothing steps on all levels.
We take enough levels to make the coarsest problem sufficiently small.
At the coarsest level, we do not solve the subproblem exactly, instead we simply just run $N_s$ iterations of proximal gradient update.
In the test we take $N_s$ iterations of (accelerated) proximal gradient update for all the pre-smoothing and post-smoothing step.
On the $\tau$ vector, we do no tuning, on all the levels we simply just take $0$ for the subdifferential in \eqref{eqn:subdiff_prox_maxpen}.
For the line search of the coarse stepsize, we simply run naive line search (Algorithm~\ref{aglo:naivelinesearch}, see Section~\ref{sec:mgprox:subsec:alpha_tune} for the discussion)
 up to machine accuracy $10^{-16}$.

\paragraph{Benchmark methods}
We compare \MGProx with the followings.

\noindent
\fbox{\texttt{ProxGrad}}: the standard proximal gradient update 
$x_0^{k+1} = \prox_{\hat{L}^{-1}g_0} \Big( x_0^k -\hat{L}^{-1} \nabla f_0(x_0^k)\Big)$, where $\hat{L}$ (initialized as $\sqrt{3}n^2/h$ based on Proposition~\ref{prop:strcvx}) is an estimate of the true Lipschitz parameter $L_0$ of $\nabla f_0$ (i.e., $\grad f_0$) at $x^k$, obtained by backtracking line search \cite[page 283]{beck2017first}.
The convergence rates of \texttt{ProxGrad} is $F(x^k) - F^* \leq \cO\Big( \min \Big\{ 1/k, \exp\big( -\frac{\mu_0}{L_0} k   \big)\Big\}\Big)$.

\noindent
\fbox{\texttt{FISTA}} \cite{nesterov1983method}: it is \texttt{ProxGrad} with Nesterov's acceleration.
With an additional auxiliary sequence $\{y_0^k\}$, the update iteration is 
    $
    x_0^{k+1} = 
    \prox_{\hat{L}^{-1}g_0} \big( y_0^k -  \hat{L}^{-1} \nabla f_0(y_0^k) \big)
    $   
    together with extrapolation
    $
    y_0^{k+1} = x_0^{k+1} + \beta_k \big( x_0^{k+1} - x_0^{k} \big)
    $.
    The stepsize $\alpha^k$ is selected using backtracking line search \cite[page 291]{beck2017first} and the extrapolation parameter $\beta_k$ follows the standard FISTA setup.
   It is proved \cite{beck2009fast} that \texttt{FISTA} has the the optimal first-order convergence rate \cite{nesterov1983method} 
   $F(x^k) - F^* \leq \cO\Big( \min \Big\{ 1/k^2 , \exp\big(-\sqrt{\frac{\mu_0}{L_0}}k \big)  \Big\} \Big)$.

\noindent
\fbox{{\texttt{FISTA-r}}}
Restarts can improve the convergence of FISTA \cite{o2015adaptive}.
There are various restarting schemes: Necoara's restart \cite{necoara2019linear}, 
adaptive restart \cite{lin2014adaptive}, and parameter-free restart \cite{aujol2023parameter}.
See \cite{aujol2021fista} for a survey.
We note that by Remark~\eqref{remark:Kmu_unknown}, both $L_0$ and $\mu_0$ of $f_0$ in \eqref{prob:EOP-1} are unknown, hence schemes \cite{necoara2019linear} and \cite{aujol2021fista} do not work.
While there are parameter-free (without knowing $L_0$ and $\mu_0$) methods \cite{lin2014adaptive}, \cite{aujol2023parameter}.
However these methods potentially requires several evaluations of $F_0$ and $\partial F$, with a cost of $\cO(N^2)$ for \eqref{prob:EOP-1}.
In the test, we implemented a simple function value restarting FISTA \cite{o2015adaptive} that only requires one evaluation of $F_0$ at every iteration, and thus is inexpensive to implement.
The convergence of FISTA-r sheds the light on the convergence of other restarting FISTA since these methods have a similar convergence rates \cite[Table 1]{aujol2021fista}.

Furthermore, there is another reason that parameter-free restarting FISTA such as \cite{lin2014adaptive} and \cite{aujol2023parameter} are not applicable for EOP.
The improvement in convergence rate from knowledge of $\mu_0$ is insubstantial if $\mu_0$ is close to $0$.
And indeed this is true for \eqref{prob:EOP-1}: we have $\mu_0$ tends to 0 as the number of elements in the mesh increases.

\noindent
\fbox{\texttt{Kocvara3}} \cite[Algorithm 3]{kocvara2016first}: a multigrid method with convergence guarantee but no convergence rate.
It is a FAS version of Kornhuber's truncation \cite{kornhuber1994monotone}, which is originally a Newton-MG (see \cref{sec:intro:subsec:review}).
 We pick the ``non-truncated version'' as it is shown to have  a better convergence.
 We remark that the method is designed for box-constrained quadratic program, and does not apply directly to \eqref{prob:EOP-1}.
        Thus, in our implementation, we adapt the structure of the algorithm (which is similar to that of \texttt{MGProx}), but with a small change: instead of projection step we perform a proximal step, since the proximal operator generalizes the projection operator.
Also, we remark that such implementation of \texttt{Kocvara3} is a special case of \MGProx that the tau vector 
$\tau$ is defined as \eqref{def:tau_mgprox2} with the set-valued subdifferentials $\partial g_0, \partial g_1$ all set to zero.

\begin{figure}[ht]
\centering
\includegraphics[width=0.99\textwidth]{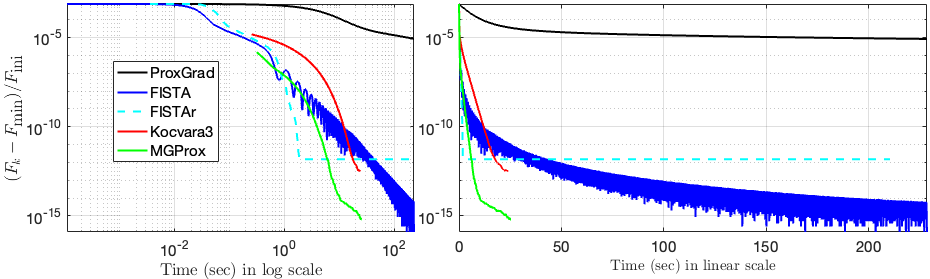}
\caption{The convergence pattern of the algorithms for the case $N^2=225$.
Left: x-axis in log scale.
Right: x-axis in linear scale.
FISTAr has the fastest convergence in the first 10 seconds, then MGProx has the fastest convergence over all.
}
\label{fig:EOP}
\end{figure}

\paragraph{The results}
Table \ref{table:EOP} shows the results on three problems with different number of variables, and Fig.\ref{fig:EOP} shows the typical convergence curves of the algorithms.
The x-axis in the plot is in terms of time and not the number of iterations, since each method has a different per-iteration cost due to the uncertainty in the number of iterations taken in the backtracking line search.

\begin{itemize}
\item In the three tests, \MGProx has the fastest convergence.
We remark that the test problem is strongly convex but the parameter $\mu_0$ is unknown and not used in any tuning in the \MGProx algorithm.
This shows case that \MGProx empirically works for convex but not strongly convex problem.

\item In general \MGProx has a higher per-iteration cost than other methods, see the discussion in Section~\ref{sec:multilevel:cost}.
However, note that the per-iteration cost of \texttt{ProxGrad} and \texttt{FISTA} is possibly expensive due to the backtracking line search.

\item We remark that we have a similar conclusion about the convergence (i.e., \MGProx has the best convergence performance in the test) if we use a larger penalty parameter in the experiment, or we if change the problem from the $\ell_1$-penalized form \eqref{prob:EOP-1} to the box-constrained form and/or the surface tension cost function is replaced by its Taylor's approximation. 

\item How the multi-level process enhances convergence speed remains an open problem.
\begin{itemize}
    \item We remark that the convergence rate of \MGProx is a loose bound since we only showed a descent condition but not a stronger sufficient descent condition.

    \item An educated guess is that the multigrid process acts as a variable metric method or pre-conditioner that reduces the number of iterations needed.
    For the EOP test problem, the evaluation of function and gradients is expensive, hence the effect of reducing number of iterations by the multi-grid process out-weights the heavier per-iteration cost of the multi-grid process.
\end{itemize}
\item As \MGProx already has the best performance among the algorithms, we do not implement \texttt{FastMGProx}.
\end{itemize}

\begin{table}[ht]\label{table:EOP}
\caption{Convergence results for \eqref{prob:EOP-1}.
The regularization parameter is $10^{-7}$ for $N^2 = 225$ and $10^{-6}$ for $N^2=65025$.
}
\centering
\begin{tabular}{ll|lll}
$N^2$ & \multicolumn{1}{c|}{Method} & \multicolumn{1}{c}{iterations $k$} & \multicolumn{1}{c}{time (sec.) } & \multicolumn{1}{c}{$\big( F_0(x_0^k) - F^{\textit{min}}_0 \Big) / F_0(x_0^{\text{ini}})$} 
\\
\hline
\parbox[t]{2mm}{\multirow{5}{*}{\rotatebox[origin=c]{90}{$(2^4-1)^2 = 225$}}} 
&\texttt{ProxGrad}           & $>10^7$ &  $232$ & $8.14 \times 10^{-6}$
\\
&\texttt{FISTA}              & $9.67 \times 10^6$  & $228$ & $1.26 \times 10^{-16}$
\\
&\texttt{FISTA-r}            & $7.77 \times 10^5$  & $210$ & $1.48 \times 10^{-12}$
\\
&\texttt{Kocvara3} $N_s = 3000$& $71$ & $24$ & $3.30 \times 10^{-13}$
\\ 
&\texttt{MGProx} $N_s = 3000$  &$50$ & $25$  & $6.31 \times 10^{-16}$
\\~\\
\hline 
\parbox[t]{2mm}{\multirow{5}{*}{\rotatebox[origin=c]{90}{$(2^8-1)^2 = 65025$}}} 
&\texttt{ProxGrad} & $>10^5$  & $335.9$ & $8.33 \times 10^{-8}$
\\
&\texttt{FISTA} &  $>10^5$  & $332.62$ & $6.64 \times 10^{-8}$
\\
&\texttt{FISTA-r} &  $>10^5$  & $364.89$ & $6.64 \times 10^{-8}$ 
\\
&\texttt{Kocvara3} $N_s = 100$&  $>10^3$ & $986.53$ & $8.07 \times 10^{-8}$
\\ 
&\texttt{MGProx} $N_s = 100$ & $50$ & $48.37$  & $1.32 \times 10^{-10}$
\\~\\
\end{tabular}
\end{table}

\section{Conclusion}\label{sec:conc}
In this work we study the combination of proximal gradient descent and multigrid method for solving a class of possibly non-smooth strongly convex optimization problems.
We propose the \MGProx method, introduce the adaptive restriction operator and provide theoretical convergence results.
We also combine \MGProx with Nesterov's acceleration, together with the optimal convergence rate with respect to the first-order methods in the function-gradient model.
Numerical results confirm the efficiency of \MGProx compared with other methods for solving Elastic Obstacle Problems.

\paragraph{Future works}
There are several problems remaining open.
\begin{itemize}
    \item The convergence rate of \MGProx is probably not tight; we conjectture that there is a tighter bound.

    \item How the multi-level process enhances the convergence speed remains open.

    \item The assumption of strong convexity of \MGProx may be relaxed.

    \item On the compactness of subdifferential, it is interesting to generalize Lemma~\ref{lemma:alphaexists} for $g_0 : \IRn \rightarrow \overline{\IR}$.
    See Remark~\ref{remark:compact}.

    \item Efficient tuning strategy for the $\tau$ selection in the subdifferential could be developed.
\end{itemize}

\section*{Acknowledgments}
We thank the three referees for their helpful comments.

\section{Appendix}\label{sec:appendix}
\subsection{Convergence of MGProx}
We present the proofs of the convergence of \texttt{MGProx}.

\subsubsection{The proof of Lemma~\ref{lemma:mgprox_suff_descent}}
%
%
\begin{proof}
By convexity and $L$-smoothness of $f$, for all $y^{k+1}, x^k, \xi$ we have
\[
\begin{array}{rlr}
f(y^{k+1})
&\leq
f(x^k) + \langle \nabla f(x^k) , y^{k+1}-x^k\rangle   + \frac{L}{2} \| y^{k+1} - x^k \|_2^2
& f \text{ is }L\text{-smooth} \dots (i)
\\
f(x^k)
&\leq
f(\xi) - \langle \nabla f(x^k) , \xi - x^k \rangle  
& f \text{ is convex} \dots (ii)
\\
f(y^{k+1}) 
&\leq
f(\xi) - \langle \nabla f(x^k) , \xi- y^{k+1} \rangle   + \frac{L}{2} \| y^{k+1} - x^k \|_2^2
& (i) + (ii)
\\
&=
f(\xi) - \big\langle \nabla f(x^k) , \xi - x^k + \frac{1}{L}G(x^k) \big\rangle   + \frac{1}{2L} \| G(x^k) \|_2^2
& y^{k+1} = x^k - \frac{1}{L}G(x^k)
\end{array}
\]
where $G(x^k)$ is the proximal gradient map of $F$ at $x^k$, see \eqref{eqn:suffdescentPGD}.

Next, adding $g(y^{k+1}) = g\big( x^k - \frac{1}{L}G(x^k) \big)$ on the both sides of the last inequality gives
\begin{equation}
F(y^{k+1}) 
\leq
f(\xi) 
- \Big\langle \nabla f(x^k) , \xi - x^k + \frac{1}{L}G(x^k) \Big\rangle 
+ \frac{1}{2L} \| G(x^k) \|_2^2
+ g\Big( x^k - \frac{1}{L}G(x^k) \Big).
\label{eqn:proof:zyx}
\end{equation}
Based on the  properties of the coarse correction (Theorem~\ref{thm:coarse_dir_descent} and Lemma~\ref{lemma:alphaexists})
and the sufficient descent property of the proximal gradient update \eqref{eqn:suffdescentPGD}, we have 
\[
F(x^{k+1}) 
~\overset{\eqref{eqn:suffdescentPGD}}{\leq}~
F(z^{k+1}) 
~\overset{\cref{thm:coarse_dir_descent}, \,\cref{lemma:alphaexists}}{\leq}~
F(y^{k+1}),
\]
so we can replace $F(y^{k+1}) $ in \eqref{eqn:proof:zyx} by $F(x^{k+1})$ and obtain
\begin{equation}
F(x^{k+1}) 
\leq
f(\xi) 
- \Big\langle \nabla f(x^k), \,  \xi - x^k + \frac{1}{L}G(x^k) \Big\rangle 
+ \frac{1}{2L} \| G(x^k) \|_2^2
+ g\Big( x^k - \frac{1}{L}G(x^k) \Big).
\label{lemma:mgprox_suff_descent1}
\end{equation}
In the following we deal with the term $g\Big( x^k - \frac{1}{L}G(x^k) \Big)$ in \eqref{lemma:mgprox_suff_descent1}.
First, by the convexity of $g$,  for all $\xi$ we have
\[
\begin{array}{lrcl}
& g(\xi) 
&\geq& \displaystyle
g\Big( x^k - \frac{1}{L}G(x^k) \Big)
+ \Big\langle\,  \partial g \big(x^k - \frac{1}{L}G(x^k) \big),\, \xi - \big(x^k - \frac{1}{L}G(x^k) \big) 
\,\Big\rangle 
\vspace{0.1mm}
\\
\iff & 
g\Big( x^k - \frac{1}{L}G(x^k) \Big)  
&\leq& \displaystyle
g(\xi) 
- \Big\langle\, \partial g \big(x^k - \frac{1}{L}G(x^k) \big),\, \xi - \big(x^k - \frac{1}{L}G(x^k) \big) \,\Big\rangle.
\end{array}
\]
By the subgradient optimality of the proximal subproblem associated with $g$, we can show that $G(x^k) - \nabla f(x^k) \in \partial g \big(x^k - \frac{1}{L}G(x^k) \big) $, hence
\begin{equation}
g\Big( x^k - \frac{1}{L}G(x^k) \Big)  
~~\leq~~
g(\xi) 
- \Big\langle G(x^k) - \nabla f(x^k),  \,  
\xi -\big(x^k - \frac{1}{L}G(x^k) \big) \Big\rangle.
\label{lemma:mgprox_suff_descent2}
\end{equation}
Combining \eqref{lemma:mgprox_suff_descent1} and \eqref{lemma:mgprox_suff_descent2} with $\xi=x^* \coloneqq \argmin F$ gives
\[
\begin{array}{rcl}
F(x^{k+1}) 
&\leq& \displaystyle
F^* 
- \Big\langle G(x^k) ,\, x^* - x^k + \frac{1}{L}G(x^k) \Big\rangle
+ \frac{1}{2L} \| G(x^k) \|_2^2
\vspace{0.1mm}
\\
&=&\displaystyle
F^* 
- \Big\langle G(x^k) ,\,  x^* - x^k \Big\rangle 
- \frac{1}{2L} \| G(x^k) \|_2^2
\vspace{0.1mm}
\\
&=& \displaystyle
F^* 
+ 
\frac{L}{2}\Big(
\| x^k - x^* \|_2^2 - \| x^k - x^* - \frac{1}{L}G(x^k) \|_2^2
\Big) 
\quad\overset{x^k - \frac{1}{L}G(x^k) \,\eqqcolon\, y^{k+1}}{\iff}~~~
\eqref{lemma:mgprox_suff_descent_iq}
\end{array}
\]
where completing the squares is used in the second equal sign.
\end{proof}

\subsection{The proof of Lemma~\ref{lemma:quadOverestimator}}
%
%
\begin{proof}
By the convexity of $f$ and $g$,
\begin{align}
f(x)\geq&\,
f(x^k) + \langle \nabla f(x^k) , x - x^k \rangle  
& f \text{ is convex} \dots (i)
\nonumber
\\
g(x) \geq&\,
g(y^{k+1}) + \langle \partial g (y^{k+1}) , x- y^{k+1} \rangle  
& g \text{ is convex} \dots (ii)
\nonumber
\\
F(x)\geq
&\,
f(x^k) + \langle \nabla f(x^k) , x - x^k \rangle  
+ g(y^{k+1}) + \langle \partial g (y^{k+1}) , x- y^{k+1} \rangle
& (i)+(ii)
\label{lemma:QI1}
\end{align}
By definitions \eqref{updt:proxgrad}, \eqref{def:prox_opertor}, the proximal gradient iteration is a majorization-minimization process that updates $x^k$ based on minimizing a local quadratic overestimator $Q$ of $x^k$, i.e., $y^{k+1} = \prox_{\frac{1}{L}g}\big( x^k - \frac{1}{L}\nabla f(x^k)\big)$ is equivalent to 
\begin{equation}
y^{k+1} 
~=~ 
\underset{\xi }{\argmin} \, 
\Big\{~
Q(\xi ; x^k) 
~\coloneqq~
f(x^k) + \big\langle \nabla f(x^k) , \xi  - x^k \big\rangle
+ \dfrac{L}{2} \|\xi  - x^k \|_2^2 
+ g(\xi )
~\Big\}.
\label{lemma:IQ_MM}
\end{equation}
Being an overestimator (which comes from the majorization-minimization interpretation \cite[Section 4.2.1]{parikh2014proximal}), we have $F(x) \leq Q(x;x^k)$, which implies for all $x$
\begin{equation}
\begin{array}{cl}
& \hspace{-2mm} F(x) - F(y^{k+1}) 
\\
\geq& \hspace{-2mm}
F(x) - Q(y^{k+1} ; x^k)
\\
\overset{\eqref{lemma:IQ_MM}}{=}& \hspace{-2mm}
F(x) 
- f(x^k) 
- \big\langle \nabla f(x^k) , y^{k+1}- x^k \big\rangle
- \dfrac{L}{2} \| y^{k+1} - x^k \|_2^2 
- g(y^{k+1})
\\
\overset{\eqref{lemma:QI1}}{\geq}& \hspace{-2mm}
\langle \nabla f(x^k) , x - x^k \rangle  
+\langle \partial g (y^{k+1}) , x- y^{k+1} \rangle
- \big\langle \nabla f(x^k) , y^{k+1} - x^k \big\rangle
- \dfrac{L}{2} \| y^{k+1} - x^k \|_2^2 
\\
=& \hspace{-2mm}
\langle \nabla f(x^k) + \partial g (y^{k+1}),  x - y^{k+1} \rangle  
- \dfrac{L}{2}\| x^k - y^{k+1} \|_2^2.
\end{array}
\label{lemma:IQ2}
\end{equation}
Applying the subgradient optimality condition to \eqref{lemma:IQ_MM} at $y^{k+1}$ gives
\[
0 ~\in~ \nabla f(x^k) + L(y^{k+1}-x^k) + \partial g (y^{k+1})
~~~\iff~~~
L(x^k-y^{k+1}) ~\in~ \nabla f(x^k)  + \partial g (y^{k+1}),
\]
so $L(x^k-y^{k+1})$ can be substituted in the first term of 
 the last line of \eqref{lemma:IQ2} and we have
\[
\begin{array}{rcl}
F(x) - F(y^{k+1})
&\geq&
L
\langle x^k-y^{k+1},  x - y^{k+1} \rangle  
- \dfrac{L}{2}\| x^k - y^{k+1} \|_2^2
\\
&=&
L\langle x^k-y^{k+1},  x - x^k + x^k - y^{k+1} \rangle  
- \dfrac{L}{2}\| x^k - y^{k+1} \|_2^2
\quad \iff \eqref{lemma:QI}.
\end{array}
\]
\end{proof}

\subsection{Convergence of FastMGProx}
We make use of Nesterov's estimate sequence to prove the convergence rate on the minimization problem
\begin{equation}\label{min:fg_appendix}
    \argmin_{x} \Big\{F(x) \coloneqq f(x) + g(x)\Big\},
\end{equation}
where $f=f_0$ and $g=g_0$ as defined in \eqref{prob:minfg}.
Since this subsection is long, to ease notation we dropped all the subscript.
Everything in this subsection refers to the fine problem.

\paragraph{Some terminology and Nesterov's estimate sequence}
\begin{itemize}
\item A quadratic over-estimator model of $f$ at $y^k$
\begin{equation}\label{eqn:f_model}
m_k(x; y^k)
~\coloneqq~
f(y^k) + \big\langle \nabla f(y^k), x - y^k \big\rangle  + \dfrac{1}{2L}\| x - y^k\|_2^2.
\end{equation}

\item A nonsmooth over-estimator model of $F$ at $y^k$
\begin{equation}\label{eqn:F_Model}
M_k(x; y^k)
~\coloneqq~
f(y^k) + \big\langle \nabla f(y^k), x - y^k \big\rangle  + \dfrac{1}{2L}\| x - y^k\|_2^2 + g(x).
\end{equation}

\item Minimizer of $M_k$ at $y^k$
\begin{equation}\label{eqn:xstar_argminM}
x^*_k(y^k)
~\coloneqq~
\argmin_x M_k(x; y^k)
\end{equation}
with the subgradient 1st-order optimality condition as
\[
0 
~\in~
\partial M_k(x^*_k(y^k); y^k) = \nabla f(y^k) +  \dfrac{x^*_k(y^k) - y^k}{L} + \partial g(x^*_k(y^k)),
\]
which implies a ``gradient'' from $x^*_k(y^k)$ and $y^k$ can be defined as
\begin{equation}\label{eqn:g}
g^k
~\coloneqq~
\dfrac{y^k - x^*_k(y^k)}{L} ~\in~ \nabla f(y^k) + \partial g(x^*_k(y^k))
\end{equation}

\item \textbf{Definition (Nesterov's estimate sequence)} \cite[Def. 2.2.1]{nesterov2003introductory}
A pair of sequences $\{\phi^k(x), \lambda^k\}$ is an \textit{estimate sequence} of an function $f$ if for all $k \in \IN$ we have
\begin{subequations}
\begin{equation}\label{Def0}
\lambda^k ~\geq~ 0,
\tag{Def 0}
\end{equation}
\begin{equation}\label{Def1}
\lambda^k ~\rightarrow~ 0,
\tag{Def 1}
\end{equation}
\begin{equation}\label{Def2}
\phi^k(x) ~\leq~ 
(1-\lambda^k)f(x) + \lambda^k \phi^0(x)
,
\tag{Def 2}
\end{equation}
\end{subequations}
where $\phi^0$ is free and $\lambda^k$ denotes $\lambda$ at the $k$th iteration, not $\lambda$ to the $k$th power.
\end{itemize}

%
%
\paragraph{How to construct an estimate sequence}
The following lemma is analogous to \cite[Lemma 2.2.2]{nesterov2003introductory}, summarizing how to construct an estimate sequence for $F$ in \eqref{min:fg_appendix}.
Compared with \cite[Lemma 2.2.2]{nesterov2003introductory} which was proposed for smooth convex optimization, the new thing here is the introduction of \eqref{A0} and the modification in \eqref{A7}.
%
%
\begin{lemma}\label{lem:how_to_construct_ES}
Assuming
\begin{subequations}
\begin{equation}\label{A0}
F(x^*_k(y^k))
~\leq~
M_k\big(x^*_k(y^k);y^k\big).
\tag{A0}
\end{equation}
\begin{equation}\label{A1}
f \text{ is }L\text{-smooth and convex (possibly not strongly convex)}.
\tag{A1}
\end{equation}
\begin{equation}\label{A2}
\phi^0(x) \text{ is a convex function}.
\tag{A2}
\end{equation}
\begin{equation}\label{A3}
\big\{ y^k \big\} \text{ is an arbitrary sequence}.
\tag{A3}
\end{equation}
\begin{equation}\label{A4a}
\big\{ \alpha^k \big\} \text{ is a sequence that } \alpha^k \in \,\,]\,0, 1\,[.
\tag{A4a}
\end{equation}
\begin{equation}\label{A4b}
\big\{ \alpha^k \big\} \text{ is a sequence that } \sum_{k=0}^\infty \alpha^k = \infty.
\tag{A4b}
\end{equation}
\begin{equation}\label{A5}
\lambda^0
~\coloneqq~
1.
\tag{A5}
\end{equation}
\begin{equation}\label{A6}
\lambda^{k+1} 
~\coloneqq~
(1-\alpha^k)\lambda^k .
\tag{A6}
\end{equation}
\begin{equation}\label{A7}
\phi^{k+1}(x) 
~\coloneqq~
(1-\alpha^k)\phi^k(x) 
+ \alpha^k \Big[
    F\big(x^*_k(y^k)\big) 
    + \big\langle g^k, x - y^k \big\rangle 
    + \frac{1}{2L}\| g^k\|_2^2
\Big] .
\tag{A7}
\end{equation}
\end{subequations}
Then the sequences $\{\phi^k(x), \lambda^k\}$ generated as in \eqref{A6}, \eqref{A7} is an estimate sequence of $F$.
\begin{proof}
The proof has 3 parts: showing $\{\lambda^k\}_{k \in \IN}$ satisfies \eqref{Def0} and \eqref{Def1}, and showing $\{\phi^k(x)\}_{k \in \IN}$ satisfies \eqref{Def2}.

\paragraph{(Showing $\lambda^k > 0$)}
We have $\lambda^{k+1} \overset{\eqref{A6}}{=} \displaystyle \prod_{k=0}^\infty (1-\alpha^k) \lambda^0\overset{\eqref{A5}}{=} \displaystyle \prod_{k=0}^\infty (1-\alpha^k) \overset{\eqref{A4a}}{>} 0$ and we showed the sequence $\{\lambda^k\}_{k \in \IN}$ is lower bounded.

\paragraph{(Showing $\lambda^k \rightarrow 0$)}
By $\lambda^{k+1} \overset{\eqref{A6}}{=} (1 - \alpha^k) \lambda^k$ we have 
$\dfrac{\lambda^{k+1}}{\lambda^k} = 1 - \alpha^k \overset{\eqref{A4a}}{<} 1$,
so $\{\lambda^k\}_{k \in \IN}$ is a monotonic decreasing sequence.
By the monotone convergence theorem, we have $\{\lambda^k\}_{k \in \IN}$ converges to $\inf \{\lambda^k\}_{k \in \IN} = 0$.
I.e., $\lambda^k \rightarrow 0$.

\paragraph{(Showing $\{\phi^k(x)\}_{k \in \IN}$ satisfies \eqref{Def2})}
We prove by induction.
The base case is true by definition and \eqref{A5}.
Assume the induction hypothesis $\phi^k(x) ~\leq~ 
(1-\lambda^k)F(x) + \lambda^k \phi^0(x)$.
Now for the case $k+1$:
\[
\begin{array}{rcll}
\phi^{k+1}(x) 
&\coloneqq&
\displaystyle
(1-\alpha^k)\phi^k(x) 
+ \alpha^k \Big[
    F\big(x^*_k(y^k)\big) 
    + \big\langle g^k, x - y^k \big\rangle 
    + \frac{1}{2L}\| g^k\|_2^2
\Big] 
& \text{by }(A7)
\vspace{0.5mm}
\\
&\leq&
\displaystyle
(1-\alpha^k)\phi^k(x) 
+ \alpha^k F(x)
& \text{($\ast$)}
\vspace{0.5mm}
\\
&=&
\displaystyle
(1-\alpha^k)\Big(\phi^k(x) + (1-\lambda^k)F(x) - (1-\lambda^k)F(x)\Big)
+ \alpha^k F(x)
\vspace{0.5mm}
\\
&=&
\displaystyle
(1-\alpha^k)\Big(\phi^k(x) - (1-\lambda^k)F(x)\Big)
+ \Big((1-\alpha^k)(1-\lambda^k)+\alpha^k\Big) F(x)
\vspace{0.5mm}
\\
&\leq&
\displaystyle
(1-\alpha^k)\Big(\lambda^k \phi^0(x) \Big)
+ \Big(1-(1-\alpha^k)\lambda^k\Big) F(x)
& \text{by induction}
\vspace{0.5mm}
\\
&=&
\displaystyle
(1-\alpha^k)\lambda^k \phi^0(x) 
+ (1-\lambda^{k+1}) F(x)
& \text{by }\eqref{A6}
\vspace{0.5mm}
\\
&<&
\displaystyle
\lambda^k \phi^0(x) 
+ (1-\lambda^{k+1}) F(x)
& \text{by }\eqref{A4a}
\end{array}
\]
where $(\ast)$ is true because we have 
\[
F(x) 
~\geq~
F\big(x^*_y(y^k)\big) 
+ L\big\langle y^k - x^*_k(y^k) , x - y^* \big\rangle 
+ \dfrac{L}{2} \| y^k - x^*_k(y^k)\|_2^2
\]
which is basically \eqref{lemma:mgprox_suff_descent1} in Lemma~\ref{lemma:mgprox_suff_descent}.
\end{proof}
\end{lemma}

As Nesterov stated in \cite{nesterov2003introductory}, up to here we are free to select the function $\phi^0(x)$.
The next lemma is about canonical closed-form expression of $\phi^{k+1}(x)$, which resembles \cite[Lemma 6]{karimi2017imro} that generalized \cite[Lemma 2.2.3]{nesterov2003introductory}.
%
%
\begin{lemma}\label{lem:canonical_phi}
Let $\phi^0(x) \coloneqq F(x^0) + \dfrac{\gamma^0}{2}\| x - z^0\|_2^2$.
Then $\phi^{k+1}$ generated recursively as \eqref{A7} in Lemma \ref{lem:how_to_construct_ES} has a closed-form expression
\begin{equation}\label{eqn:phi_k1}
    \phi^{k+1}(x) ~=~ \overline{\phi}^{k+1} + \dfrac{\gamma^{k+1}}{2} \| x - z^{k+1} \|_2^2,
\end{equation}
where 
\begin{subequations}
\begin{equation}\label{eqn:lem:canonical_phi_gamma}
\gamma^{k+1} = (1-\alpha^k) \gamma^k,
\end{equation}
\begin{equation}\label{eqn:lem:canonical_phi_z}
z^{k+1} = z^k - \dfrac{\alpha^k}{\gamma^{k+1}} g^k,
\end{equation}
\begin{equation}\label{eqn:lem:canonical_phi_bar}
\overline{\phi}^{k+1} = (1-\alpha^k)\overline{\phi}^k
+ \alpha^k F\big(
x^*_k(y^k)
\big)
+ 
\dfrac{\alpha^k}{2}
\Big(
\dfrac{1}{L}
- \dfrac{\alpha^k}{\gamma^{k+1}}
\Big)
\| g^k\|_2^2
+ 
\alpha^k \big\langle g^k, z^k - y^k \big\rangle.
\end{equation}
\end{subequations}
\begin{proof}
A 4-part proof.   
\paragraph{Part 1. Proving the sequence $\{\phi^k(x)\}$ satisfies \eqref{eqn:phi_k1}}
We first prove $\nabla^2 \phi^k(x) = \gamma^k I$ for all $k$ by induction.
The base case at $k=0$ holds by taking the Hessian of $\phi^0(x) \coloneqq F(x^0) + \frac{\gamma^0}{2}\| x - z^0\|_2^2$ from the assumption.
Assuming the induction hypothesis $\nabla^2 \phi^{k}(x) = \gamma^{k} I$, now for the case $k+1$, taking the Hessian of $\phi^{k+1}$ in \eqref{A7} gives
\begin{equation}\label{lem:canonical_phi_quadratic_form_intermediate}
\nabla^2 \phi^{k+1}(x)
~\overset{\eqref{A7}}{\coloneqq}~
(1-\alpha^k) \nabla^2  \phi^k(x) 
\overset{\text{induction hypothesis}}{=}
(1-\alpha^k) \gamma^k I
~\overset{\eqref{eqn:lem:canonical_phi_gamma}}{=}~
\gamma^{k+1} I.
\end{equation}
Eq.\eqref{lem:canonical_phi_quadratic_form_intermediate} implies that for all $k$ the function $\phi^k$ is a quadratic function in the form of \eqref{eqn:phi_k1}, for some constants $\overline{\phi}^k, \gamma^k$, $z^k$.

\paragraph{Part 2. Proving the sequence $\{\gamma^k\}$ satisfies \eqref{eqn:lem:canonical_phi_gamma}}
This is proved by the last equality in \eqref{lem:canonical_phi_quadratic_form_intermediate}.

\paragraph{Part 3. Proving the sequence $\{z^k\}$ satisfies \eqref{eqn:lem:canonical_phi_z}}
First we consider the gradient and Hessian of $\phi^{k+1}(x)$.
\begin{subequations}
\begin{equation}\label{eqn:A7_phi_k1}
\phi^{k+1}(x)
~\overset{\eqref{A7},\eqref{eqn:phi_k1}}{\coloneqq}~
(1-\alpha^k)\Big[  \overline{\phi}^k + \dfrac{\gamma^k}{2} \| x - z^k \|_2^2  \Big]
+ \alpha^k \Big[    F\big(x^*_k(y^k)\big)     + \big\langle g^k, x - y^k \big\rangle     + \frac{1}{2L}\| g^k\|_2^2 \Big], 
\end{equation}
\begin{equation}\label{eqn:grad_A7_phi_k1}
\nabla \phi^{k+1}(x)
~\overset{\eqref{eqn:A7_phi_k1}}{=}~
(1-\alpha^k) \gamma^k (x-z^k) + \alpha^k  g^k.
\end{equation}
\end{subequations}
Setting $\nabla \phi^{k+1}(x) = 0 $ in \eqref{eqn:grad_A7_phi_k1} gives an expression of $z^k$
\begin{equation}\label{lem:canonical_phi_z_zero}
\nabla \phi^{k+1}(x) ~\overset{\eqref{eqn:grad_A7_phi_k1}}{=}~ (1-\alpha^k) \gamma^k (x-z^k) + \alpha^k  g^k ~=~ 0.
\end{equation}
Now we take the gradient of $\phi^{k+1}$ using \eqref{eqn:phi_k1} gives
\begin{equation}\label{lem:canonical_phi_quadratic_z_zero_k1}
\nabla \phi^{k+1}(x) ~=~ \gamma^{k+1}(x - z^{k+1}) = 0.
\end{equation}
Combine \eqref{lem:canonical_phi_z_zero}, \eqref{lem:canonical_phi_quadratic_z_zero_k1} and use \eqref{eqn:lem:canonical_phi_gamma} gives \eqref{eqn:lem:canonical_phi_z}.

\paragraph{Part 4. Proving the sequence $\{\overline{\phi}^{k+1}\}$ satisfies \eqref{eqn:lem:canonical_phi_bar}}
Equating \eqref{eqn:phi_k1} and \eqref{eqn:A7_phi_k1} at $x=y^k$ gives
\begin{equation}\label{eqn:lem:canonical_phi_last2}
\overline{\phi}^{k+1} + \dfrac{\gamma^{k+1}}{2} \| y^k - z^{k+1} \|_2^2
\overset{\eqref{eqn:phi_k1}}{=}
\phi^{k+1}(y^k) 
\overset{\eqref{eqn:A7_phi_k1}}{=}
(1-\alpha^k)\Big[  \overline{\phi}^k + \dfrac{\gamma^k}{2} \| y^k - z^k \|_2^2  \Big]
+ \alpha^k \Big[    F\big(x^*_k(y^k)\big)     + \frac{1}{2L}\| g^k\|_2^2 \Big].
\end{equation}
Note that 
\begin{equation}\label{eqn:lem:canonical_phi_last1}
\dfrac{\gamma^{k+1}}{2} \| z^{k+1} -  y^k \|_2^2
~\overset{\eqref{eqn:lem:canonical_phi_z}}{=}~
\dfrac{\gamma^{k+1}}{2} \| z^k -  y^k \|_2^2
- \gamma^{k+1} \big\langle z^k - y^k, \dfrac{\alpha^k}{\gamma^{k+1}}g^k \big\rangle 
+ \dfrac{\gamma^{k+1}}{2} \| \dfrac{\alpha^k g^k}{\gamma^{k+1}} \|_2^2.
\end{equation}
Combine \eqref{eqn:lem:canonical_phi_last2}, \eqref{eqn:lem:canonical_phi_last1} and using \eqref{eqn:lem:canonical_phi_gamma} gives \eqref{eqn:lem:canonical_phi_bar}.
\end{proof}
\end{lemma}
At this stage we recall that by \eqref{A3} in Lemma~\ref{lem:how_to_construct_ES}, the sequence $\{y^k\}_{k \in \IN}$ is ``free''.
In view of Lemma~\ref{lem:canonical_phi}, we would like to construct $y^k$ such that $\overline{\phi}^{k+1} \geq F(x^*_k(y^k)) \eqqcolon F(x^{k+1})$, which holds for $k=0$ that $\overline{\phi}^0 =F(x^0)$.
That is, we define $x^{k+1} = x^*_k(y^k) \overset{\eqref{eqn:xstar_argminM}}{\coloneqq} \argmin M_k(x;y^k)$

%
%
The following lemma resembles \cite[Theorem 3]{karimi2017imro} that reveals the importance of estimate sequence and also the condition $\overline{\phi}^{k+1} \geq F(x^*_k(y^k)) \eqqcolon F(x^{k+1})$.
The lemma differs from \cite[Lemma 2.2.1]{nesterov2003introductory} on smooth convex optimization.
\begin{lemma}\label{lem:rate_bound}
For minimization problem \eqref{min:fg_appendix}, assume $x^* \in X^* \coloneqq \argmin F(x)$ exists and denote $F^* \coloneqq F(x^*)$.
Suppose $\displaystyle F(x^k) \leq \overline{\phi}^k \coloneqq \min_{x} \phi_k(x)$ holds for a sequence $\{x^k\}_{k \in \IN}$, where $\{\phi^k, \lambda^k\}_{k \in \IN} $ is an estimate sequence of $F$, and we define $\phi^0 \coloneqq F(x^0) + \dfrac{\gamma^0}{2}\| x^0 - x^* \|_2^2 $, then we have for all $k \in \IN$ that
\[
F(x^k) - F^*
~\leq~
\lambda^k \Big[ F(x^0) + \dfrac{\gamma^0}{2}\| x^0 - x^* \|_2^2 -F^*
\Big].
\]
\begin{proof} Starting from the assumption $\displaystyle F(x^k) \leq \overline{\phi}^k \coloneqq \min_{x} \phi_k(x)$,
\[
F(x^k) 
\leq
\overline{\phi}^k
~\eqqcolon~ 
\displaystyle \min_{x} \phi^k
\overset{\eqref{Def2}}{\leq}
\displaystyle \min_{x}\,
(1-\lambda^k)F(x) + \lambda^k \phi^0(x)
~\leq~
(1-\lambda^k)F^* + \lambda^k \phi^0(x^*),
\]
\[
\begin{array}{rrcl}
\implies
& F(x^k) - F^* 
&\leq& 
\lambda^k \Big[  \phi^0(x^*) - F^* \Big]
~=~ 
\lambda^k \Big[ F(x^0) + \dfrac{\gamma^0}{2}\| x^0 - x^* \|_2^2 -F^* \Big].
\end{array}
\]
\end{proof}
\end{lemma}
Lemma~\ref{lem:rate_bound} tells that the convergence rate of $\big\{F(x^k) - F^*\big\}_{k \in \IN}$ follows the convergence rate of $\big\{\lambda^k \big\}_{k \in \IN}$ with $\big\{\lambda^k \big\}_{k \in \IN} \overset{\eqref{Def0}}{>} 0$ and $\big\{\lambda^k \big\}_{k \in \IN} \overset{\eqref{Def1}}{\rightarrow} 0$, and thus finding the convergence rate of $\big\{\lambda^k \big\}_{k \in \IN}$ gives the convergence rate of $\big\{F(x^k) - F^*\big\}_{k \in \IN}$.

%
%
The following lemma resembles \cite[Lemma 2.2.1]{nesterov2003introductory}.
\begin{theorem}\label{thm:lambda_rate}
Suppose $\displaystyle F(x^k) \leq \overline{\phi}^k \coloneqq \min_{x} \phi_k(x)$ holds for a sequence $\{x^k\}_{k \in \IN}$, where $\{\phi^k, \lambda^k\}_{k \in \IN} $ is an estimate sequence of $F$.
Define $\phi^0 \coloneqq F(x^0) + \dfrac{\gamma^0}{2}\| x^0 - x^* \|_2^2 $.
Assuming all the conditions in Lemma~\ref{lem:how_to_construct_ES}, Lemma~\ref{lem:canonical_phi} and Lemma~\ref{lem:rate_bound}.
Then we have
\[
0 ~<~
\lambda^k 
~<~
\dfrac{4L}{
(2\sqrt{L} - \sqrt{\gamma^0})^2
+ 2\Big(
2\sqrt{L}\sqrt{\gamma^0}
-\sqrt{\gamma^0}
\Big) \sqrt{\gamma^0}k
+ (\sqrt{\gamma^0} k)^2
}.
\]
\begin{proof}
A long, highly-involved tedious mechanical proof.
We start with $\gamma^{k+1}$ in \eqref{eqn:lem:canonical_phi_gamma}:
\begin{equation}\label{eqn:thm:lambda_rate:gamma_k1}
\gamma^{k+1} 
~\overset{\eqref{eqn:lem:canonical_phi_gamma}}{=}~
(1-\alpha^k)\gamma^k
~\overset{\eqref{eqn:lem:canonical_phi_gamma}}{=}~
\cdots
~\overset{\eqref{eqn:lem:canonical_phi_gamma}}{=}~
\prod_{i=0}^k (1-\alpha^i)\gamma^0
~\overset{\eqref{A5},\eqref{A6}}{=}~
\lambda^{k+1}\gamma^0
.
\end{equation}
Next, by the first step in Algorithm FMGProx that 
$L(\alpha^k)^2 = (1-\alpha^k)\gamma^k \overset{\eqref{eqn:lem:canonical_phi_gamma}}{\eqqcolon} \gamma^{k+1}$, solving for the root $\alpha^k $ gives
\begin{equation}\label{eqn:thm:lambda_rate:alphak_root}
\alpha^k 
~=~ 
\pm \sqrt{\dfrac{\gamma^{k+1}}{L}}
~\overset{\eqref{eqn:thm:lambda_rate:gamma_k1}}{=}~
\pm \sqrt{\dfrac{\lambda^{k+1} \gamma^0}{L}}
.
\end{equation}
Now consider $\lambda^k$ in \eqref{A6}:
\begin{equation}\label{eqn:thm:lambda_rate:lambda}
1 - \dfrac{\lambda^{k+1}}{\lambda^k} 
~\overset{\eqref{A6}}{=}~
\alpha^k 
~\overset{\eqref{eqn:thm:lambda_rate:alphak_root}}{=}~
\pm \sqrt{\dfrac{\gamma^{k+1}}{L}}
~\overset{\eqref{eqn:thm:lambda_rate:gamma_k1}}{=}~
\pm \sqrt{\dfrac{\lambda^{k+1}\gamma^0}{L}}
.
\end{equation}
Now by \eqref{A4a},\eqref{A6} and \eqref{A7} we have for all $k$ that
$\{\lambda^k\}_{k \in \IN}$ is a \textit{strictly} nonnegative sequence
and a \textit{strictly} monotonic decreasing sequence, i.e.:
\begin{subequations}
\begin{equation}\label{eqn:thm:lambda_rate:lambda_positive}
\lambda^k \overset{\eqref{A4a},\eqref{A6}}{=} \displaystyle \prod_{i=0}^k (1 - \alpha^k) \overset{\eqref{A4a}}{>} 0
,
\end{equation}
\begin{equation}\label{eqn:thm:lambda_rate:lambda_decrease}
\lambda^{k+1} \overset{\eqref{A4a},\eqref{A6}}{<} \lambda^k
~~\overset{\eqref{eqn:thm:lambda_rate:lambda_positive}}{\iff}~~
\dfrac{1}{\sqrt{\lambda^k}} ~<~ \dfrac{1}{\sqrt{\lambda^{k+1}}}
.
\end{equation}
\end{subequations}
Hence we can divide \eqref{eqn:thm:lambda_rate:lambda} by $\lambda^{k+1} \overset{\eqref{eqn:thm:lambda_rate:lambda_positive}}{>} 0$:
\begin{equation}\label{eqn:thm:lambda_rate:lambda_diff}
\dfrac{1}{\lambda^{k+1}}  - \dfrac{1}{\lambda^k} 
~=~
\pm \dfrac{1}{\sqrt{\lambda^{k+1}}}\sqrt{\dfrac{\gamma^0}{L}}
.
\end{equation}
On the left hand side of \eqref{eqn:thm:lambda_rate:lambda_diff}, by the fact that $\lambda^k \overset{\eqref{eqn:thm:lambda_rate:lambda_positive}}{>} 0$ for all $k$, we have
\begin{equation}\label{eqn:thm:lambda_rate:lambda_diff_sqrt}
\dfrac{1}{\lambda^{k+1}}  - \dfrac{1}{\lambda^k} 
~=~
\Big(
\dfrac{1}{\sqrt{\lambda^{k+1}}}  + \dfrac{1}{\sqrt{\lambda^k}}
\Big)
\Big(
\dfrac{1}{\sqrt{\lambda^{k+1}}}  - \dfrac{1}{\sqrt{\lambda^k}}
\Big)
~\overset{\eqref{eqn:thm:lambda_rate:lambda_decrease}}{<}~
\dfrac{2}{\sqrt{\lambda^{k+1}}}  
\Big(
\dfrac{1}{\sqrt{\lambda^{k+1}}}  - \dfrac{1}{\sqrt{\lambda^k}}
\Big)
.
\end{equation}
Now combine \eqref{eqn:thm:lambda_rate:lambda_diff} and \eqref{eqn:thm:lambda_rate:lambda_diff_sqrt} with $\gamma^0 > 0$ and $L > 0$ gives 
\[
\displaystyle \pm \dfrac{1}{\sqrt{\lambda^{k+1}}}\sqrt{\dfrac{\gamma^0}{L}}
<
\dfrac{2}{\sqrt{\lambda^{k+1}}}  
\Big(
\dfrac{1}{\sqrt{\lambda^{k+1}}}  - \dfrac{1}{\sqrt{\lambda^k}}
\Big)
~\implies~
- \dfrac{1}{2} \sqrt{\dfrac{\gamma^0}{L}}
<
\dfrac{1}{2} \sqrt{\dfrac{\gamma^0}{L}}
<
\dfrac{1}{\sqrt{\lambda^{k+1}}}
- \dfrac{1}{\sqrt{\lambda^{k}}}.
\]
Let $\theta^k \coloneqq \dfrac{1}{\sqrt{\lambda^k}}$ be an increasing sequence (since $\lambda^k$ is decreasing), we have
\[
\theta^{k+1}
~>~
\theta^{k} +
\dfrac{1}{2} \sqrt{\dfrac{\gamma^0}{L}}
~>~
\theta^{k-1} +
\dfrac{1}{2} \sqrt{\dfrac{\gamma^0}{L}} + \dfrac{1}{2} \sqrt{\dfrac{\gamma^0}{L}}
~>~ \cdots 
~\overset{\lambda^0 = 1}{>}~
1 + \dfrac{k}{2}\sqrt{\dfrac{\gamma^0}{L}}.
\]
Thus $
\dfrac{1}{\sqrt{\lambda^{k+1}}}= \theta^{k+1}
>
1 + \dfrac{k}{2}\sqrt{\dfrac{\gamma^0}{L}}
$ gives
$
0 ~\overset{\eqref{Def0}}{<}~\sqrt{\lambda^{k+1}} ~<~ \dfrac{1}{1 + \dfrac{k}{2}\sqrt{\dfrac{\gamma^0}{L}}}
$ and hence
\[
0 
~<~ 
\lambda^{k+1} 
~<~
\dfrac{1}{1 + k \sqrt{\dfrac{\gamma^0}{L}} + \dfrac{k^2}{4} \dfrac{\gamma^0}{L} }
~~\implies~~
0 
<
\lambda^k
<
\dfrac{1}{1 + (k-1)\sqrt{\dfrac{\gamma^0}{L}} + \dfrac{(k-1)^2}{4}\dfrac{\gamma^0}{L}}
\]
\[
~~\implies~~
0 
<
\lambda^k
<
\dfrac{4L}{
(2\sqrt{L} - \sqrt{\gamma^0})^2
+ 2\Big(
2\sqrt{L}\sqrt{\gamma^0}
-\sqrt{\gamma^0}
\Big) \sqrt{\gamma^0}k
+ (\sqrt{\gamma^0} k)^2
}.
\]
\end{proof}
\end{theorem}

Lastly, the following corollary gives the proof of Theorem~\ref{thm:FMGProx_fastr_ate} in the main text.
\begin{corollary}\label{corollary:FMGProx_final_convergence_rate}
For the sequence $\{x^k\}$ produced by Algorithm~FMGProx, we have 
\[
\begin{array}{rcl}
F(x^k) - F^*  
~~\leq~~
\dfrac{4L
\Big(
F(x^0) + \dfrac{\gamma^0}{2}\| x^0 - x^* \|_2^2 - F^*
\Big)
}{
(2\sqrt{L} - \sqrt{\gamma^0})^2
+ 2\gamma^0\Big(2\sqrt{L}- 1\Big) k
+ (\sqrt{\gamma^0} k)^2
}
.
\end{array}
\]
\begin{proof}
Combine Lemma~\ref{lem:rate_bound} and Theorem~\ref{thm:lambda_rate}.
\end{proof}
\end{corollary}

\bibliographystyle{siamplain}
\bibliography{references}
\end{document}